\documentclass[10pt,a4paper,twoside,bibliography=totocnumbered]{article}

\usepackage{xcolor}

\usepackage{graphicx}

\usepackage{amsmath}
\usepackage{amssymb}

\usepackage{geometry}

\usepackage{pdflscape}

\begin{document}

\title{{\huge \textcolor{black}{A Full-Density Approach to Simulating Random Iteration Equations with Applications}\vspace*{-0.4cm}\\ \begin{minipage}{5cm}\centering \textcolor{black}{---}\vspace*{-0.4cm}
\end{minipage}\\}{\Large \textcolor{black}{Original Research Article}}\vspace*{0.5cm}}


\author{Wolfgang Hoegele$^{1}$\vspace*{0.4cm}\\
{\normalsize $^{1}$ Munich University of Applied Sciences HM}\\{\normalsize Department of Computer Science and Mathematics}\\{\normalsize Lothstraße 64, 80335 München, Germany}\vspace*{0.4cm}\\
{\normalsize corresponding mail: \texttt{wolfgang.hoegele@hm.edu}}\vspace*{0.4cm}\\
ORCID: 0000-0002-5303-9334 \vspace*{0.4cm}}

\date{\today}

\maketitle
\thispagestyle{empty}

\textit{This manuscript is published and can be cited by}\medskip

\textit{W. Hoegele, A Full-Density Approach to Simulating Random Iteration Equations with Applications,\\ Communications in Nonlinear Science and Numerical Simulation, 163, part 1. \\ https://doi.org/10.1016/j.cnsns.2026.110508.}\medskip

\textit{This manuscript appears also on ArXiv.org\\
arXiv:2603.17466 [math.DS],\\ https://doi.org/10.48550/arXiv.2603.17466}

\section*{Abstract}

The goal of this study is to introduce a unified framework for simulating random iteration equations (RIE), understood as iteration equations containing random variables. {The main idea is to propagate approximations of the full state density from one iteration to the next, rather than estimating it from many repeated pathwise Monte Carlo simulations.} The presentation of the {\textit{RIE modeling framework} is conceptually simple} based on recent work on static random equations and designed to be accessible. The modeling requirements for RIEs allow for potential nonsmooth nonlinearities and stochasticities in the transfer function. {Additionally, the \textit{RIE computational strategy} for full-density propagation is presented based on iterative likelihood / posterior calculations.} As results, illustrative applications of nonlinear random and stochastic differential equation simulations, a new full-density gradient descent method (FDGD) for global optimization under uncertainty and examples of chaotic mappings are presented in order to demonstrate the breadth of the utility of this framework. In total, the character of the presentation is explorative and encourages new applications and theoretical studies. \medskip

\textbf{Keywords:} random iteration equation, density propagation, uncertainty quantification, nonlinear dynamics, full-density gradient descent, stochastic differential equation \medskip

\newpage
\thispagestyle{empty}
\tableofcontents

\section*{About the Author}

Dr. Högele is Professor of Applied Mathematics and Computational Science at the Department of Computer Science and Mathematics at the Munich University of Applied Sciences (HM), Germany. His research interests are in general mathematical modeling and, more specifically, stochastic modeling, simulation, and analysis of complex systems in applied mathematics.\medskip

\thispagestyle{empty}

\newpage

\section{Introduction}

In this work a new {modeling and} computational framework is presented to investigate the evolution of \textit{random iteration equations} (RIE). In general, first order iteration equations are transforming a state vector $\boldsymbol{x}\in\mathbb{R}^R$ starting from an initial state $\boldsymbol{x}^{(1)}$ with a stepwise iteration mapping from the $n$-th iteration state $\boldsymbol{x}^{(n)}$ to the $(n+1)$-th iteration state $\boldsymbol{x}^{(n+1)}$. This is denoted by the transfer function $\boldsymbol{T} : \mathbb{R}^R \rightarrow \mathbb{R}^R$ with $\boldsymbol{x}^{(n+1)} = \boldsymbol{T}(\boldsymbol{x}^{(n)})$. RIEs will be defined as iteration equations that contain \textit{random variables}, i.e. $\boldsymbol{T}$ contains stochasticities summarized by random variables in the iteration equation and, in consequence, the resulting next state $\boldsymbol{x}^{(n+1)}$ must also be interpreted as a random variable with a probability density. This means, the RIE takes the density of the current state and computes the density of the next state by incorporating the random variable densities in the transfer function. In total, we are observing the \textit{stepwise evolution of state probability density functions} from iteration to iteration, being an example of a \textit{density propagation} algorithm \cite{Li 2006,Caluya 2020}. The proposed framework shows how these density functions can be computed in practice which is demonstrated by expressive examples. RIEs can be regarded as a key component for the simulation of many dynamical systems that transform density functions, which places this study in the broad field of \textit{uncertainty quantification} \cite{Sullivan2015,Garcia-Sanchez2022}. A minor point in this argumentation, but a major difference in the application, is that also the initial state $\boldsymbol{x}^{(1)}$ will be treated as a random variable, i.e. the iteration begins not with single starting values, but with a full distribution. {Standard approaches to simulate dynamics numerically depending on uncertainties or diffusion processes are based on statistical pathwise Monte Carlo simulation, i.e. from all probability densities a single sample is drawn and the corresponding deterministic iteration is performed. This pathwise approach is performed in repetition and all results are collected statistically in a histogram map in order to see the final density of paths \cite{Rubinstein2017,Kloeden1992}. This differs from direct density propagation: pathwise Monte Carlo estimates densities statistically from sampled trajectories rather than propagating densities directly.}

{This work has two primary objectives: a) Introducing the \textit{RIE modeling framework}: Demonstrating that working with RIEs as the central argument is a unifying flexible modeling framework to connect random/stochastic dynamics of very different application fields. b) Presenting the \textit{RIE computational strategy}: Proposing an approximate Monte Carlo integration method for propagation of the full state-space probability densities under random iteration maps using a smoothed transition kernel. Specifically, this method focuses on the transient phase of the stochastic simulation (possibly including nonsmooth nonlinearities and stochasticities in the transfer function) by computing a regularized transfer.}

{With respect to a), in order to highlight the broader relevance of the modeling framework, three major application areas of iteration equations in applied mathematics are presented}: simulation of dynamical systems (including \textit{random} and \textit{stochastic differential equations}), optimization (extending \textit{gradient descent} optimization to a \textit{full-density gradient descent} method) and chaotic maps (presenting the strange attractors of the \textit{Ikeda map} and the \textit{Lozi map}). This is not the full {modeling} potential of this framework, instead it shows its relevance for important areas in applied mathematics utilizing expressive example simulations. 

{Giving context to b),} based on a recent publication which focused on a method to compute the probabilistic solution space of \textit{{static} random equations} \cite{Hoegele 2026}, the propagation of the probability densities can be performed directly by computing the full state density from iteration to iteration considering the complete stochastic environment {of the transfer function}. {The computational approach utilizes a likelihood / posterior reformulation}. Although in the numerical framework we utilize Monte Carlo integration techniques, this should not be misinterpreted with the statistical collection of pathwise Monte Carlo simulations. 

{There are general boundaries of this work, which we state explicitly: First, this work highlights the flexibility of the \textit{RIE modeling framework}, but it does not argue its general superiority compared to other perspectives on stochastic dynamics. We are presenting illustrative example applications which are not definitive for the field but indicate broad plausibility and versatility. Second, the proposed \textit{RIE computational strategy} for full-density propagation should be regarded as an alternative or complementary approach to standard stochastic simulation methods, such as pathwise Monte Carlo simulation and comparison for verification will be presented in the course of the paper. Third, this work does not focus on discretization schemes or discuss discretization errors in detail, e.g. how to get from random or stochastic differential equations to RIEs, instead it focuses on how to model and compute the RIEs as the starting point once it is agreed on a discretization scheme.  Fourth, we typically utilize the normal distribution for the parameter random variables or diffusion processes which is not a limitation of the computational framework but nevertheless not explicitly treated in the paper. Non-standard probability densities (also non-standard diffusion processes) combined with significant nonlinearities or discontinuities or without Lipschitz contraction properties of the transfer mapping remain interesting applications for the future.}

{\subsection{Related Approaches}}

Major contributions on the convergence theorems {for stationary / invariant densities} are presented in the closely related field of \textit{iterated random functions} focusing on existence statements utilizing Lipschitz contraction arguments \cite{Sandric 2025, Diaconis 1999}. {Such stationary investigations are regarded as a different objective compared to this work} since our approach is focusing on the \textit{transient phase} and explicit evolution of the density propagation through stepwise nonlinear iterations. Specifically, this study demonstrates that the density propagation of stochastic and dynamic iteration problems can be simulated straightforwardly, even if the theoretical treatment might be difficult or convergence is not guaranteed or even violated. In our perspective, this is more aligned with many real world applications in which mostly only a limited number of iterations may be of computational interest.

Specifically with respect to the theoretical treatment of stochastic dynamical systems this approach is therefore a computational complement to \textit{Perron-Frobenius} operator approaches, e.g. \cite{Surasinghe 2024}, (which are typically applied to deterministic, non-random transfer functions and focus on the limiting density, i.e. the invariant measure) and \textit{Fokker-Planck-Equation} or \textit{Kolmogorov forward equation} approaches, e.g. see \cite{Oksendal 2013}, (which are difficult to apply if continuous-time Gaussian diffusion processes or smooth drift or diffusion coefficients are not guaranteed). And on the computational side it is also complementary to \textit{pathwise Monte Carlo} simulations due to the focus on \textit{full-density propagation}. {For a supplementary exploration of the theoretical context, the detailed connection to the closely related \textit{Perron-Frobenius} operator (Section \ref{sec:GenForm}) and \textit{Ulam's} method for its approximation (Appendix \ref{sec:AdvAnaRitz}) are presented. }

{\subsection{Outline of the Paper}

The methods section is structured as follows: In Section \ref{sec:GenForm} the general description of the \textit{RIE modeling framework} as well as the reformulation as a likelihood / posterior density are derived, including the reformulation as an \textit{integral linear operator} and a theoretical example for illustration. In Section \ref{sec:numimp} the \textit{RIE computational strategy} is presented based on direct Monte Carlo integration of the posterior density, including a contextualization with an alternative, supplementary computational approach related to \textit{Ulam's method} which is further discussed in Appendix  \ref{sec:AdvAnaRitz}. Finally, two major application areas are connected with the general \textit{RIE modeling framework}: In Section \ref{sec:meth_RDESDE} the theoretical connection to the simulation of random and stochastic differential equations is presented and in Section \ref{sec:meth_FDGD} the utilization for a stochastic extension of the gradient descent algorithm in optimization is demonstrated which is named as the \textit{full-density gradient descent}.

In the results Section \ref{sec:res:dynamicsys} it is demonstrated how this modeling framework is applied to dynamical system evaluation utilizing the simulation of predator prey \textit{Rosenzweig McArthur} random and stochastic differential equations. As a supplementary verification in Appendix \ref{sec:NumVec2DOrnUhl} the analytic solution of the \textit{2D Ornstein-Uhlenbeck SDE} is compared to the presented RIE computation as well as to pathwise Monte Carlo simulation. In Section \ref{res:FDGD} the stochastically extended gradient descent is applied to a saddle point objective function as well as the challenging \textit{Himmelblau's} function. Finally, in Section \ref{sec:ChaotiCRIE}  the simulation of the strange attractors of the chaotic \textit{Ikeda map} and the \textit{Lozi map} are demonstrated utilizing the RIE computation.

}

\section{Methods}

\subsection{General Formulation}
\label{sec:GenForm}
We introduce the \textit{random iteration equation} (RIE) with an at least piecewise continuous transfer function $\boldsymbol{T}: \mathbb{R}^R\times \mathbb{R}^K \rightarrow \mathbb{R}^R$ and a random variable parameter vector $\boldsymbol{C}^{(n)}$ with $K$ entries and known densities (depending on the iteration number $n$) by
\begin{align}
\boldsymbol{x}^{(n+1)} = \boldsymbol{T}( \boldsymbol{x}^{(n)}, \boldsymbol{C}^{(n)} )\;,
\end{align}
given an initial state $\boldsymbol{x}^{(1)}\sim f_{\boldsymbol{x}^{(1)}}$. We want to utilize the methodology recently introduced in \cite{Hoegele 2026} to rewrite this into the standard \textit{random equation} form $\boldsymbol{M}(\boldsymbol{x};\boldsymbol{A}) = \boldsymbol{B}$ with independent random variable vectors $\boldsymbol{A}$ and $\boldsymbol{B}$, and $\boldsymbol{x}$ being the corresponding \textit{solution vector}. We transform the iteration equation to
\begin{align}
\boldsymbol{x}^{(n+1)} - \boldsymbol{T}( \boldsymbol{x}^{(n)}, \boldsymbol{C}^{(n)} ) = \boldsymbol{0}
\end{align}
and by comparison this leads to $\boldsymbol{B}$ a random variable vector centered around zero and with small standard deviations compared to other standard deviations, the solution vector $\boldsymbol{x}^{(n+1)}$ for each iteration and the two random variables $\boldsymbol{x}^{(n)}$ and $\boldsymbol{C}^{(n)}$. What a sufficiently small standard deviation for the zero vector $\boldsymbol{B}$ is, needs to be investigated for each application and can be better understood in the following theoretical example. As a minimal requirement piecewise continuity of $\boldsymbol{T}$ is stated, which is regarded as a comparably weak assumption and sufficient for many applications explicitly including nonsmooth transfer functions. It should be noted that even weaker measure theoretic requirements can be stated which are inherited directly from the integrability assumptions about $\boldsymbol{M}$ presented in \cite{Hoegele 2026} and are based on \textit{Lebesgue measurability}. In total, this leads to
\begin{align}
&\boldsymbol{M} : \mathbb{R}^R\times \mathbb{R}^{R+K} \rightarrow \mathbb{R}^R\\
&\boldsymbol{M}(\boldsymbol{x}^{(n+1)};\underbrace{\boldsymbol{x}^{(n)}, \boldsymbol{C}^{(n)}}_{=:\boldsymbol{A}}) := \boldsymbol{x}^{(n+1)} - \boldsymbol{T}( \boldsymbol{x}^{(n)}, \boldsymbol{C}^{(n)} )
\end{align}
with $\boldsymbol{A}:=(x_1^{(n)},\dots,x_R^{(n)}, C_1^{(n)}, \dots, C_{K}^{(n)})^T$ the random variable vector consisting of two parts where $\boldsymbol{x}^{(n)}$ is a $R$-dimensional random vector with density $f_{\boldsymbol{x}^{(n)}}$ computed in the $n$-th iteration, and $\boldsymbol{C}^{(n)}$ a static $K$-dimensional random vector with \textit{i.i.d.} $C_1^{(1)},\dots, C_1^{(n)}\sim f_{C_1}$ for the first entry, until \textit{i.i.d} $ C_{K}^{(1)},\dots,C_{K}^{(n)}\sim f_{C_{K}}$ for the last entry. The common density function is denoted by $f_{\boldsymbol{C}}$. For a well-defined iteration a starting density $f_{\boldsymbol{x}^{(1)}}$ must be defined. The independence of the $C_i^{(n)}$ for different $n$ is directly in the spirit of deterministic iteration equations, where the transfer function itself is independent of previous evaluations of it.

Since $\boldsymbol{x}^{(n)}$ and $\boldsymbol{C}^{(n)}$ are independent, we get a likelihood distribution for $\boldsymbol{x}^{(n+1)}$ based on those densities and the specific transfer function $\boldsymbol{T}$ by
\begin{align}
\mathcal{L}_{\boldsymbol{x}^{(n+1)}}(\boldsymbol{0}|\boldsymbol{x}) = \int\limits_{\mathbb{R}^{K}}  \int\limits_{\mathbb{R}^R} f_{\boldsymbol{B}}(  \boldsymbol{x} - \boldsymbol{T}( \boldsymbol{s}_1, \boldsymbol{s}_2 ) ) \cdot f_{\boldsymbol{x}^{(n)}}(\boldsymbol{s}_1)\cdot f_{\boldsymbol{C}}(\boldsymbol{s}_2)\;\text{d}\boldsymbol{s}_1\,\text{d}\boldsymbol{s}_2\;.\label{equ:GenIntegral}
\end{align}
whose general form was presented previously \cite{Hoegele 2026}. {This formula points out the importance of the regularization by $\boldsymbol{B}$ for the computational strategy.} When focusing on probability densities, we reformulate the likelihood function with Bayes' theorem to a posterior density by utilizing a non- or weakly-informative prior in this study which effectively scales the likelihood intensities 
\begin{align}
\pi_{\boldsymbol{x}^{(n+1)}}(\boldsymbol{x})\propto \mathcal{L}_{\boldsymbol{x}^{(n+1)}}(\boldsymbol{0}|\boldsymbol{x})\;,
\end{align}
essentially normalizing the likelihood function to an integral value of $1$. In consequence, we approximate the next-step density at each iteration via the normalized posterior $f_{\boldsymbol{x}^{(n+1)}} := \pi_{\boldsymbol{x}^{(n+1)}}$, evolving from a prior starting density $f_{\boldsymbol{x}^{(1)}}:=\pi_{\boldsymbol{x}^{(1)}}$. Of course, in each iteration there could also be a density modifying prior utilized which strongly influences the density propagation. This could be important, for example, if the state space is restricted by hard or soft boundaries (utilizing a uniform prior in the region of admissible solutions) and/or is a discrete space with discrete probability densities (utilizing a gridded prior). We will leave such interesting possibilities outside of this study (although it is not outside of the framework) and only refer to the corresponding Section 2.5 in \cite{Hoegele 2026}.

Rewriting Equation \ref{equ:GenIntegral} into
\begin{align}
\pi_{\boldsymbol{x}^{(n+1)}}(\boldsymbol{x})\propto  \int\limits_{\mathbb{R}^R} \underbrace{\left(\; \int\limits_{\mathbb{R}^{K}}  f_{\boldsymbol{B}}(  \boldsymbol{x} - \boldsymbol{T}( \boldsymbol{s}_1, \boldsymbol{s}_2 ) ) \cdot f_{\boldsymbol{C}}(\boldsymbol{s}_2)\,\text{d}\boldsymbol{s}_2 \right)}_{=:k(\boldsymbol{x},\boldsymbol{s}_1)}  \cdot\; \pi_{\boldsymbol{x}^{(n)}}(\boldsymbol{s}_1)\;\text{d}\boldsymbol{s}_1\;,
\label{equ:RIEIntOp}
\end{align}
the RIE iteration can be reframed as an integral linear operator acting on the space of posterior densities with kernel $k(\cdot,\cdot)$, which we denote as the \textit{linear RIE operator}. We emphasize that this reformulation is in general possible including cases with nonlinearities, discontinuities or non-standard parameter probability densities in the  transfer function $\boldsymbol{T}$. These \textit{RIE kernels} are generally bounded and positive, since they consist of probability density functions, but stronger properties are not guaranteed. Spectral investigation for specific kernels may lead to further insights into the dynamic properties of the density propagation for selected RIEs. {A direct comparison with the deterministic \textit{Perron-Frobenius} operator is presented using its integral operator representation \cite{Surasinghe 2024}:
\begin{align}
\pi_{\boldsymbol{x}^{(n+1)}}(\boldsymbol{x})\propto  \int\limits_{\mathbb{R}^R}	 \delta( \boldsymbol{x} - \boldsymbol{T}( \boldsymbol{s} ) ) \cdot\; \pi_{\boldsymbol{x}^{(n)}}(\boldsymbol{s})\;\text{d}\boldsymbol{s}\;,
\end{align}
where $\delta$ denotes the delta distribution. This comparison emphasizes that the \textit{linear RIE operator} can be interpreted as a regularized / smoothed stochastic extension of the \textit{Perron-Frobenius} operator, or shorter, a \textit{regularized stochastic transfer operator}, incorporating general stochasticity in the transfer (not only additive noise as typically used \cite{Surasinghe 2024}). Regularization via the zero random variable $\boldsymbol{B}$ is necessary for formulating the likelihood and in the numerical simulation as discussed in Section \ref{sec:numimp}.}

\subsection*{Theoretical Example}

In order to gain intuition about this operator, we consider a simple and artificial iteration, which reproduces its one-dimensional initial density of $x^{(1)}$ by adding diffusion with the zero mean random variables $C^{(n)}$
\begin{align}
x^{(n+1)} := x^{(n)} + C^{(n)}\;\Rightarrow\;T(x^{(n)},C^{(n)}) := x^{(n)} + C^{(n)}\;,
\end{align}
which leads to the \textit{linear RIE operator}
\begin{align}
\pi_{x^{(n+1)}}(x)\propto\;&  \int\limits_{\mathbb{R}} \left(\; \int\limits_{\mathbb{R}}  f_{B}(  x - s_1 - s_2 ) ) \cdot f_{C}(s_2)\,\text{d} s_2 \right)  \cdot\; \pi_{x^{(n)}}(s_1)\;\text{d} s_1\;\\
=\;&  \int\limits_{\mathbb{R}} f_{B\ast C}(  x - s_1 )  \cdot\; \pi_{x^{(n)}}(s_1)\;\text{d} s_1\\
=\;&  f_{B\ast C\ast x^{(n)}}(  x  ) \propto f_{(B\ast C)^n \ast x^{(1)}}(  x  )\;.
\end{align}
This is obviously an iterative convolution of the initial density with the dominant diffusion density and the minor zero value density $B$, and exactly what we intuitively expect considering the iteration equation.

\subsection{Numerical Implementation}
\label{sec:numimp}

{A main goal of this work is to present a framework for calculating the iterations of the RIEs.} As presented in \cite{Hoegele 2026}, the integral in Equation \ref{equ:GenIntegral}, can be approximated by Monte Carlo integration, utilizing $P$ samples with the formula
\begin{align}
\pi_{\boldsymbol{x}^{(n+1)}}(\boldsymbol{x})\propto \;  \frac{1}{P} \sum\limits_{p=1}^P\left(\prod\limits_{r=1}^R f_{B_r}( x_r - T_r( \boldsymbol{s}_{1,p}, \boldsymbol{s}_{2,p} ))\right)\;,\label{equ:NumCalcMC}
\end{align}
where the samples $\boldsymbol{s}_{1,p}$ are drawn form the density $\pi_{\boldsymbol{x}^{(n)}}$ and $\boldsymbol{s}_{2,p}$ are drawn from $f_{\boldsymbol{C}^{(n)}}$. The sampling for the latter is straightforward, since these are given by the RIE definition and one can for example utilize \textit{latin hypercube sampling} (even the same representative samples for all iterations once calculated before iterations start). The sampling of the first density $\pi_{\boldsymbol{x}^{(n)}}$ must be done in each iteration and these densities are non-trivial.  There are standard \textit{Monte Carlo Markov Chain }(MCMC) methods in order to draw from such posteriors, but in this work we prefer due to simplicity an \textit{acceptance rejection} approach \cite{Rubinstein2017}. Essentially, drawing in a $R$-dimensional rectangular box around the density region of  $\pi_{\boldsymbol{x}^{(n)}}(\boldsymbol{x})$ uniformly distributed samples $\boldsymbol{x}_p$ with one additional uniform sample $y_p$ in the range of $0$ and the maximum overall value of $\pi_{\boldsymbol{x}^{(n)}}$. If $y_p$ is larger than $\pi_{\boldsymbol{x}^{(n)}}(\boldsymbol{x}_p)$ it is rejected, else it is accepted. This is done in parallelization for speed-up and exactly $P$ samples are computed in each iteration. {During this study the box around the density region is kept constant and identical with the constant computational grid of the densities as they are presented in the results figures. Further, the necessary maximum density value for the acceptance rejection algorithm is taken as the maximum value of the computed density on that grid and evaluation of the density function at the randomly selected points is performed by linear interpolation.}

The major computational compromise must be found between the number of samples $P$, the standard deviations of $\boldsymbol{B}$ and the accepted noise level in the final density by utilizing Monte Carlo integration. For example, the more concentrated $\boldsymbol{B}$ is around zero (what is desirable since too large standard deviations introduce a blur on the posterior density or even slight diffusion over a large number of iterations), the more samples inserted in $f_{B_r}$ are typically not contributing to the final posterior density increasing the noise level. \medskip

{A schematic algorithm for the \textit{RIE computational strategy} (which is also the basis of the presented results)  for evaluating the smoothed full-density propagation on a calculation grid is given by}
\begin{itemize}
\item[1.] Define the RIE, especially the transfer function  $\boldsymbol{T}: \mathbb{R}^R\times \mathbb{R}^K \rightarrow \mathbb{R}^R$, the \textit{i.i.d.} random variable densities $f_{\boldsymbol{C}}$ and the initial prior density $\pi_{\boldsymbol{x}^{(1)}}$. Set $n=1$.
\item[2.] Take $P$ draws $\boldsymbol{s}_{2,p}$ ($p=1,\dots,P$) from $f_{\boldsymbol{C}}$ by standard sampling methods (e.g. by \textit{latin hypercube} sampling).
\item[3.] Take $P$ draws $\boldsymbol{s}_{1,p}$ ($p=1,\dots,P$) from $\pi_{\boldsymbol{x}^{(n)}}$ (e.g. by \textit{acceptance rejection}).
\item[4.] Evaluate Equation \ref{equ:NumCalcMC} numerically in order to calculate $\pi_{\boldsymbol{x}^{(n+1)}}$ on the same grid as $\pi_{\boldsymbol{x}^{(n)}}$ was provided (including the normalization step of the approximated density {by utilizing the Riemann sum}).
\item[5.] Set $n \rightarrow n+1$ and repeat steps 3--5 until a certain limit has been reached, such as a given iteration number or the changes between $\pi_{\boldsymbol{x}^{(n)}}$ and $\pi_{\boldsymbol{x}^{(n+1)}}$ are below a certain tolerance measure.
\end{itemize}

This algorithm for simulating full densities of the RIEs can be applied directly for the stochastic extension of a given iteration equation, such as the chaotic mappings presented in the results Section \ref{sec:ChaotiCRIE}. More important in applied mathematics are situations where the iteration equation itself (and its stochastic extension to RIE) is a tool to solve challenging (mostly nonlinear) problems. We will present two of the most prominent use cases in the next subsections.

It is pointed out that utilizing the perspective of the \textit{linear RIE operator} combined with a linear decomposition of the posterior density function as linear combination of basis functions ({a typical \textit{Galerkin} approach}), much more efficient numerical schemes are possible than the presented Monte Carlo algorithm. In this case the coefficients of the linear combination are updated in each iteration by multiplying them with a \textit{propagation matrix}. This perspective can also be utilized as an efficient method to incorporate spatial discretization directly into the algorithm by (effectively) pixel or voxel basis functions {approximating \textit{Ulam's} method. The full argumentation is presented in the Appendix \ref{sec:AdvAnaRitz}.}

\subsection{Discretization of Dynamical Systems with Stochasticity as Random Iteration Equation}
\label{sec:meth_RDESDE}

We will see, that being able to simulate RIEs is a key step to simulate approximately \textit{random differential equation} (RDE) and \textit{stochastic differential equation} (SDE) dynamical systems in order to observe their evolution of the full densities over time \cite{Jornet 2023,Kloeden1992}. {It is stated with clarity that we do not suggest underestimating the theoretically challenging investigation of RDEs and SDEs by using them only as an example application in the RIE modeling framework. On the contrary, we recognize that full depth investigation of such dynamical systems is a difficult and ongoing task of current research. Instead we just want to show, that the RIE modeling framework is flexible enough to also address the RDE and SDE simulation rational which may open up new opportunities for applications.} We are starting applications with a general ODE system with random variable parameters $\boldsymbol{C}$, representing a RDE,
\begin{align}
\dot{\boldsymbol{x}} = \boldsymbol{F}(t,\boldsymbol{x}, \boldsymbol{C})\;,
\end{align}
with a well-defined right hand side $\boldsymbol{F}$. By applying the discretization $\dot{\boldsymbol{x}} \approx \frac{\boldsymbol{x}^{(n+1)}-\boldsymbol{x}^{(n)}}{\Delta t}$, this leads to the iteration equation
\begin{align}
\boldsymbol{x}^{(n+1)} = \boldsymbol{x}^{(n)} + \Delta t\cdot \boldsymbol{F}(t_n,\boldsymbol{x}^{(n)},\boldsymbol{C}^{(n)})\;,
\end{align}
essentially representing a random \textit{Euler iteration} equation. Please note that the independent drawing of $\boldsymbol{C}^{(n)}$ in each iteration differs from the classical RDE simulation, where parameters are sampled once per path and held fixed. In this sense, the presented framework may be understood either as an approximation of a classical RDE by an RIE or as a non-classical RDE evolution in its own right. When using the nomenclature of Section \ref{sec:GenForm}, this leads to the general transfer function
\begin{align}
\boldsymbol{T}( \boldsymbol{x}^{(n)}, \boldsymbol{C}^{(n)} ) := \boldsymbol{x}^{(n)} + \Delta t\cdot \boldsymbol{F}(t_n,\boldsymbol{x}^{(n)},\boldsymbol{C}^{(n)})\;,
\end{align}
for an arbitrary dynamical system. Of course, more advanced discretization schemes can be utilized, but for simplicity we leave this demonstration applying the \textit{Euler} algorithm.\medskip

A major observation is that when utilizing this discretization scheme, we can also simulate the evolution of probability density functions of SDEs in It\^{o} convention with random coefficients   
\begin{align}
\text{d} \boldsymbol{x} = \boldsymbol{F}(t,\boldsymbol{x},\boldsymbol{G})\, \text{d} t + \boldsymbol{B}(t,\boldsymbol{x},\boldsymbol{D})\, \text{d} \boldsymbol{W}_t \;,
\end{align}
with the drift coefficient function $\boldsymbol{F}$, the diffusion coefficient function $\boldsymbol{B}$ and a Wiener process  $\boldsymbol{W}_t$, utilizing the \textit{Euler-Maruyama} algorithm \cite{Rubinstein2017,Kloeden1992,Higham 2001,JianZu2023} by introducing new terms
\begin{align}
\boldsymbol{x}^{(n+1)} = \boldsymbol{x}^{(n)} + \boldsymbol{F}(t_n,\boldsymbol{x}^{(n)},\boldsymbol{G}^{(n)})\, \Delta t + \boldsymbol{B}(t_n,\boldsymbol{x}^{(n)},\boldsymbol{D}^{(n)})\, \Delta \boldsymbol{W}^{(n)} \;,\label{equ:RSDE1}
\end{align}
with the coefficient random variable vectors $\boldsymbol{G}^{(n)}$ and $\boldsymbol{D}^{(n)}$, and the discrete Wiener process increment $\Delta \boldsymbol{W}^{(n)} =  \boldsymbol{W}_{t_{n+1}}- \boldsymbol{W}_{t_{n}}$ as normally distributed random variables with zero mean and variances $\Delta t$. This leads to the general transfer function for SDEs by renaming the random variables
\begin{align}
\boldsymbol{T}( \boldsymbol{x}^{(n)}, \boldsymbol{C}^{(n)} )  := \boldsymbol{x}^{(n)} + \boldsymbol{F}(t_n,\boldsymbol{x}^{(n)},\boldsymbol{C}_{1:K_1}^{(n)})\, \Delta t + \boldsymbol{B}(t_n,\boldsymbol{x}^{(n)},\boldsymbol{C}_{K_1+1:K_1+K_2}^{(n)})\, \boldsymbol{C}_{K_1+K_2+1:K_1+K_2+R}^{(n)}\;,\label{equ:RSDE2}
\end{align}
with three parts of $\boldsymbol{C}$, i.e. the $K_1$ random variables of the drift coefficient $\boldsymbol{C}_{1:K_1}$ $ (=\boldsymbol{G})$, the $K_2$ random variables of the diffusion coefficient $\boldsymbol{C}_{K_1+1:K_1+K_2}$  $(=\boldsymbol{D})$ and the $R$ random variables of the discrete increment in the Wiener process $\boldsymbol{C}_{K_1+K_2+1:K_1+K_2+R}$ $(=\Delta \boldsymbol{W})$.

Please note, that Equation \ref{equ:RSDE2} is a general model including parameter uncertainties (the RDE part) as well as a diffusion process (the SDE part) in one single computational framework. This shows already from the methodological side the flexibility of the RIE approach for approximating those in general challenging processes combined. Additionally, it can be observed that if $\boldsymbol{G}^{(n)}$ (and respectively $\boldsymbol{D}^{(n)}$) of Equation \ref{equ:RSDE1} are independent (but not necessarily identically distributed) for different $n$ then also more complicated coefficient processes can be simulated just by reinterpretation of the terms (not leaving the computational framework). {This exemplifies the conceptual similarities between RDEs and SDEs in this simulation perspective when applying full-density propagation. For clarification: when applying pathwise Monte Carlo simulation, RDEs and SDEs need to be simulated differently, e.g. for each path in RDEs only one fixed parameter set is drawn from the densities and deterministic calculation is performed, while for SDEs during the path calculation continuously random errors are drawn.} Further, the SDEs are not limited to Wiener processes, but can in the simulation take any other stochastic process that increment can be calculated for each iteration (which is deliberately general and leaves a lot of possibilities). Finally, also more sophisticated approaches than the simplistic \textit{Euler-Maruyama} algorithm can be utilized for SDE simulation, without leaving the presented general framework as long as they lead to iteration equations. For the convergence rates of standard SDE types for the \textit{Euler-Maruyama} algorithm see e.g. Higham\cite{Higham 2001}.

Obviously, this strict focus on the discretized iteration equation as the key quantity opens up many possibilities for simulating challenging RDEs and SDEs besides the standard theory and current research in their fields. In this study, this is only regarded as one example application of RIE simulations.

\subsection{Full-Density Gradient Descent as Random Iteration Equation}
\label{sec:meth_FDGD}

The idea of gradient descent is to update a current position $\boldsymbol{x}^{(n)}$ in a minimization process of an objective function $F(\boldsymbol{x})$ with $F:\mathbb{R}^R \rightarrow \mathbb{R}$ stepwise in the direction of the steepest downward direction which is represented by the negative direction of the local gradient, i.e. leading to the gradient descent (GD) algorithm
\begin{align}
\boldsymbol{x}^{(n+1)} = \boldsymbol{x}^{(n)} - \eta\cdot \nabla F(\boldsymbol{x}^{(n)})\;,
\end{align}
with a \textit{learning rate} $\eta$. We present several stochastic extensions to this classical optimization routine, which family will be denoted with \textit{full-density gradient descent} (FDGD):
\begin{itemize}
\item[1.] \textbf{FDGD-I:} Minimization of a stochastic objective function containing $K$ random variables for parameters of the objective function
\begin{align}
\boldsymbol{x}^{(n+1)} = \boldsymbol{x}^{(n)} - \eta\cdot \nabla F(\boldsymbol{x}^{(n)},\boldsymbol{C}^{(n)})
\end{align}
In this version \textit{stochastic optimization} is performed, i.e. we have an objective function containing random variable parameters. We just want to present this approach here since it is an obvious version, but we are not studying it in detail further on.
\item[2.] \textbf{FDGD-II:} Minimization of a deterministic objective function updating the deterministic direction simultaneously for a full density of directions, leading to $R$ random variables
\begin{align}
\boldsymbol{x}^{(n+1)} = \boldsymbol{x}^{(n)} - \eta\cdot \nabla F(\boldsymbol{x}^{(n)}) + \boldsymbol{C}^{(n)}
\end{align}
This approach introduces a \textit{diffusion} process into the iteration which helps exploring the objective function. 
\item[3.] \textbf{FDGD-III:} Additionally to FDGD-II the update utilizes a full density of learning rates, leading to $R+1$ random variables
\begin{align}
\boldsymbol{x}^{(n+1)} = \boldsymbol{x}^{(n)} - C_{R+1}^{(n)}\cdot \nabla F(\boldsymbol{x}^{(n)}) + \boldsymbol{C}_{{1:R}}^{(n)}\;,
\end{align}
leading to an update with different optimization schemes simultaneously.
\end{itemize}
FDGD-II and FDGD-III contain a diffusion component due to the added random variables in each iteration with a simultaneous convergence to local minima, producing complex dynamics during optimization. This leads, for example, to the transfer function for FDGD-III
\begin{align}
\boldsymbol{T}( \boldsymbol{x}^{(n)}, \boldsymbol{C}^{(n)} ) := \boldsymbol{x}^{(n)} - C_{R+1}^{(n)}\cdot \nabla F(\boldsymbol{x}^{(n)}) + \boldsymbol{C}_{{1:R}}^{(n)}\;,
\end{align}
and for the other FDGD schemes analogously. Again, this formulation is very general for any differentiable objective function $F$. In addition, combining FDGD-I and FDGD-III addresses an extremely broad set of stochastic optimization problems, transforming the local gradient descent algorithm to a \textit{global stochastic optimization} technique. The difference to \textit{stochastic gradient descent} (SGD) in the literature is that we do not use a single realization of the stochastic path in each iteration (which in SGD is produced by subsets of data for calculating the gradient in machine learning) but the whole density functions are transformed which allows a much broader search \cite{Goodfellow2016}. Related ideas of density propagation by stochastic gradient descent \cite{Mandt 2017} and for the connection of SDEs and gradient descent \cite{Przbylowicz 2026} are part of recent research. We further want to draw the attention to the strong similarities of SDEs and the FDGD approach within the RIE framework which may lead to further insights.

\section{Results}

\subsection{Evolution of Dynamical Systems}
\label{sec:res:dynamicsys}

\begin{figure}[htbp]
\centering
\includegraphics[width=14cm]{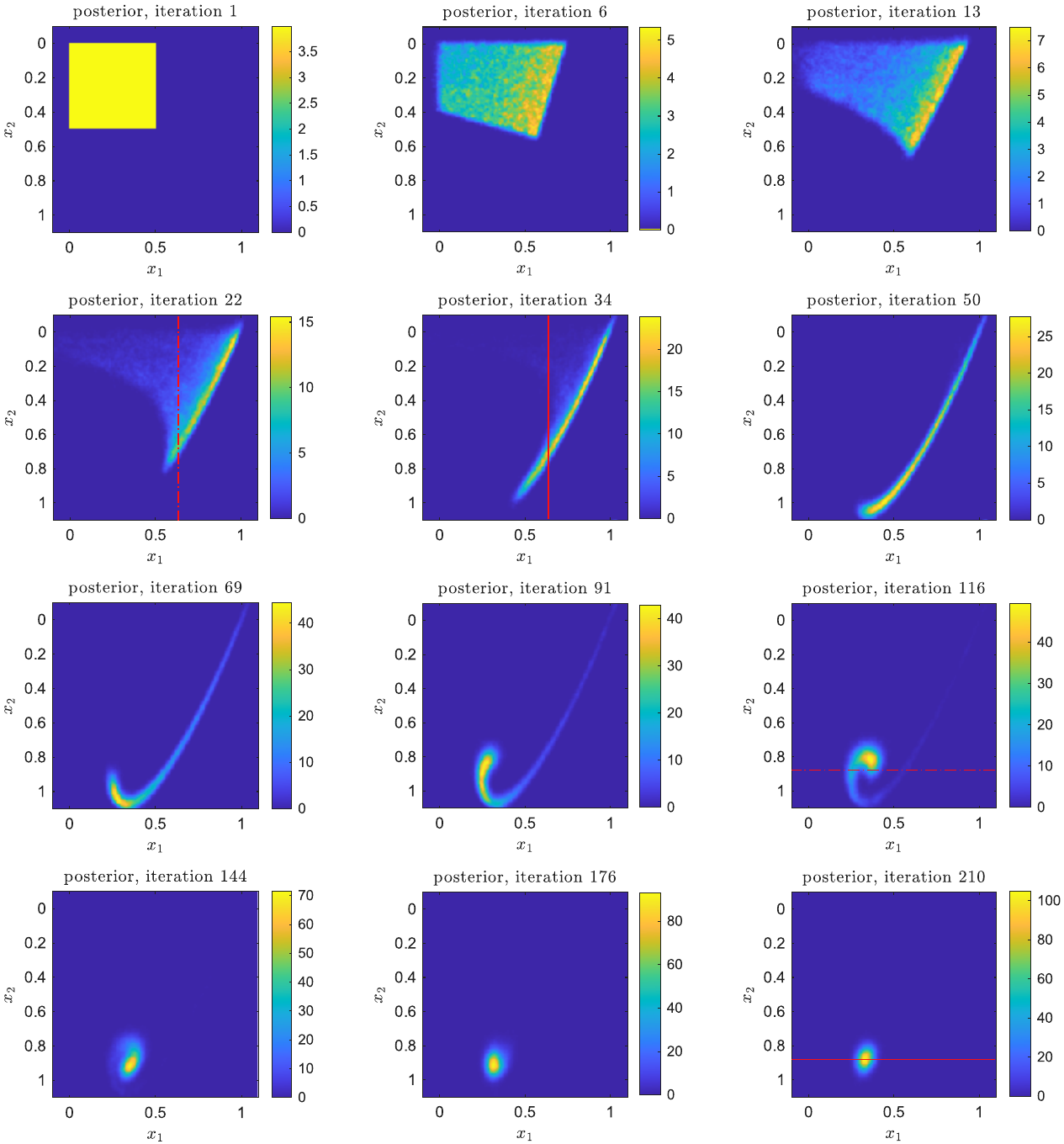}
\caption{Illustration of the posterior densities for different iteration numbers of the Rosenzweig McArthur RDE model. One can see the evolution of the square uniform initial density and how it evolves due to the dynamics of the Rosenzweig McArthur model with random variable parameters to the steady state density. {The computational density values were evaluated at $200\times 200$ grid points with equidistant spacing on $[-0.1,1.1]\times[-0.1,1.1]$. The initial density was uniform on the box $[0,0.5]\times [0,0.5]$.}}  
\label{fig:Res:RosMcA_iteration_posterior_004_2}
\end{figure}

\begin{figure}[htbp]
\centering
\includegraphics[width=14cm]{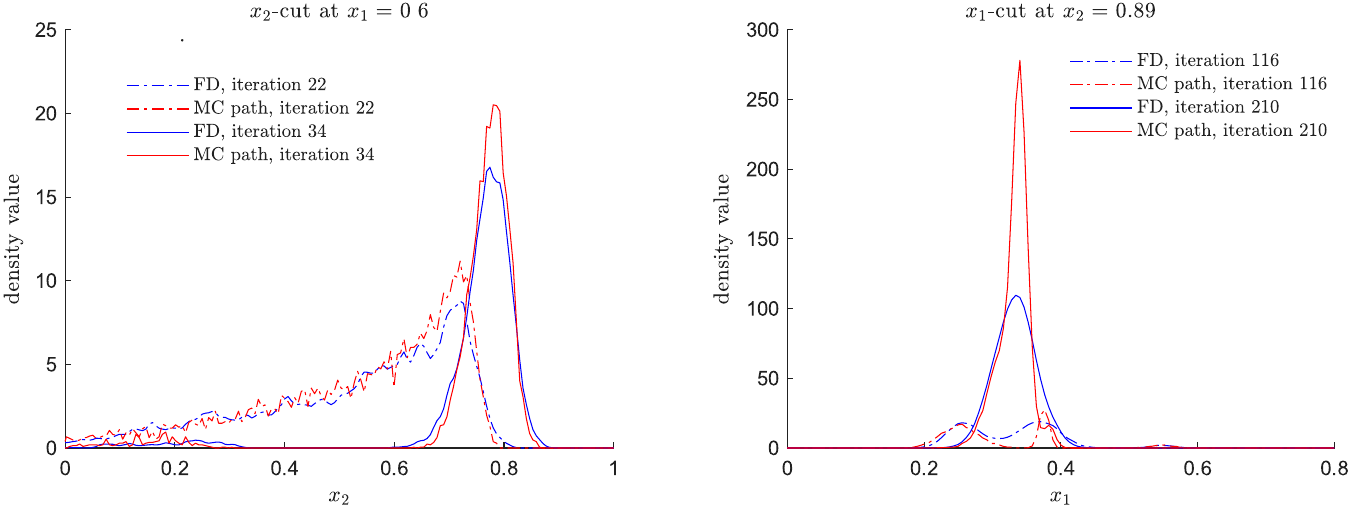}
\caption{{Verification of the computation of Figure \ref{fig:Res:RosMcA_iteration_posterior_004_2} by comparing cross sections of the numerical simulation of the full-density approach (FD) to the pathwise Monte Carlo simulation (MC path) utilizing $768000$ paths. The cross section are additionally presented in Figure \ref{fig:Res:RosMcA_iteration_posterior_004_2}. Left: $x_2$-cuts at $x_1=0.6$ for iterations $22$ and $34$, right: $x_1$-cuts at $x_2=\frac{8}{9}$ for iterations $116$ and $210$.}  }
\label{fig:Res:RosMcA_iteration_posterior_004_2_cuts}
\end{figure}

{In order to show how the RIE modeling framework and computational strategy can be applied to concrete RDE and SDE systems we utilize a highly nonlinear predator prey dynamical system, the normalized \textit{Rosenzweig McArthur} model. We present two versions of this example: First, the RDE model \cite{Stollenwerk 2022,Wyse 2022, Hoegele 2026_2} for the predator prey simulation is given by}
\begin{align}
\dot{x}_1 & = x_1\,\left(1-\frac{x_1}{C_1}\right)-\frac{C_2\,x_1\,x_2}{1+x_1}\\
\dot{x}_2 & = -C_3\,x_2 +\frac{C_2\,x_1\,x_2}{1+x_1}\;,
\end{align}
with the random variables representing uncertainties about the parameters of the model. The straightforward discretization $\dot{x} \approx \frac{x^{(n+1)}-x^{(n)}}{\Delta t}$ leads to the transfer function
\begin{align}
\boldsymbol{T}( \boldsymbol{x}^{(n)}, \boldsymbol{C}^{(n)} ) & := \left(\begin{array}{c} x_{1}^{(n)} + \Delta t\cdot \left( x_{1}^{(n)}\,\left(1-\frac{x_{1}^{(n)}}{C_1^{(n)}}\right)-\frac{C_2^{(n)}\,x_1^{(n)}\,x_2^{(n)}}{1+x_1^{(n)}} \right)\\
x_{2}^{(n)} + \Delta t\cdot \left( -C_3^{(n)}\,x_2^{(n)} +\frac{C_2^{(n)}\,x_1^{(n)}\,x_2^{(n)}}{1+x_1^{(n)}} \right) \end{array}\right)
\end{align}
In the simulation, we utilize independent Gaussian densities $f_{C_1} := \mathcal{N}(1,0.01^2)$, $f_{C_2} := \mathcal{N}(1,0.01^2)$ and $f_{C_3} := \mathcal{N}(0.25,0.01^2)$ as well as $f_{B_1},f_{B_2} := \mathcal{N}(0,0.005^2)$. {Please note, numerically the $f_{B_1},f_{B_2}$ are truncated after three times the standard deviation.} The differential equation model is known to lead to a steady state density \cite{Hoegele 2026_2} and we want to observe the evolution of the starting distribution to this steady state distribution based on the discretized iteration equation. We utilize $P=192000$ samples and $\Delta t=0.2$ (observing acceptable convergence results with minor numerical artifacts, see Figure \ref{fig:Res:RosMcA_iteration_numericstep} (left) for an impression). In Figure \ref{fig:Res:RosMcA_iteration_posterior_004_2}, we are utilizing starting density $f_{\boldsymbol{x}^{(1)}}$ as a uniform density on $[0,0.5]\times [0,0.5]$ and observe its evolution into the steady state distribution. The deterministic steady state (when taking the modes of each density function as parameter values) is at $(\frac{1}{3},\frac{8}{9})$ and coincides well with the simulation result. {In addition, in Figure \ref{fig:Res:RosMcA_iteration_posterior_004_2_cuts} cross sections of that simulation are presented and compared to pathwise Monte Carlo simulation. The shapes and positions of the intensity values are comparable but the smoothing effect of $\boldsymbol{B}$ in the full-density calculation and the RDE approximation by RIEs becomes visible.}

\begin{figure}[htbp]
\centering
\includegraphics[width=14cm]{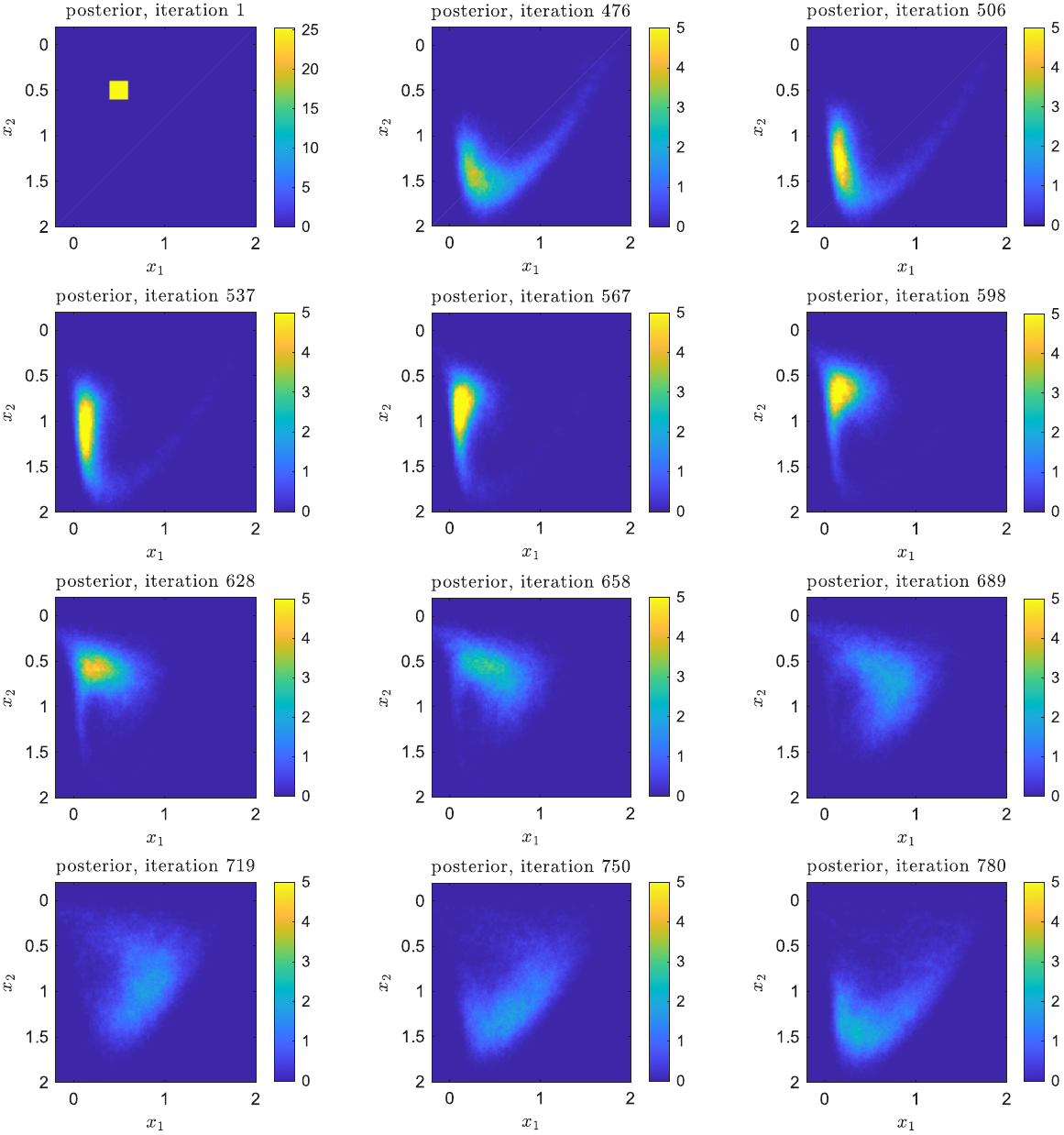}
\caption{Illustration of the posterior densities for different iteration numbers of the Rosenzweig McArthur SDE model. The evolution of the densities is presented roughly over one period of the orbit in the $x_1$-$x_2$-space starting at iteration 476 ($t=23.8$), with dominant stochastic effects due to the SDE diffusion and initial state density. {The computational density values were evaluated at $200\times 200$ grid points with equidistant spacing on $[-0.2,2]\times[-0.2,2]$. The initial density was uniform on the box $[0.4,0.6]\times [0.4,0.6]$.}}  
\label{fig:Res:RosMcA_iteration_posteriorSDE_004}
\end{figure}

\begin{figure}[htbp]
\centering
\includegraphics[width=14cm]{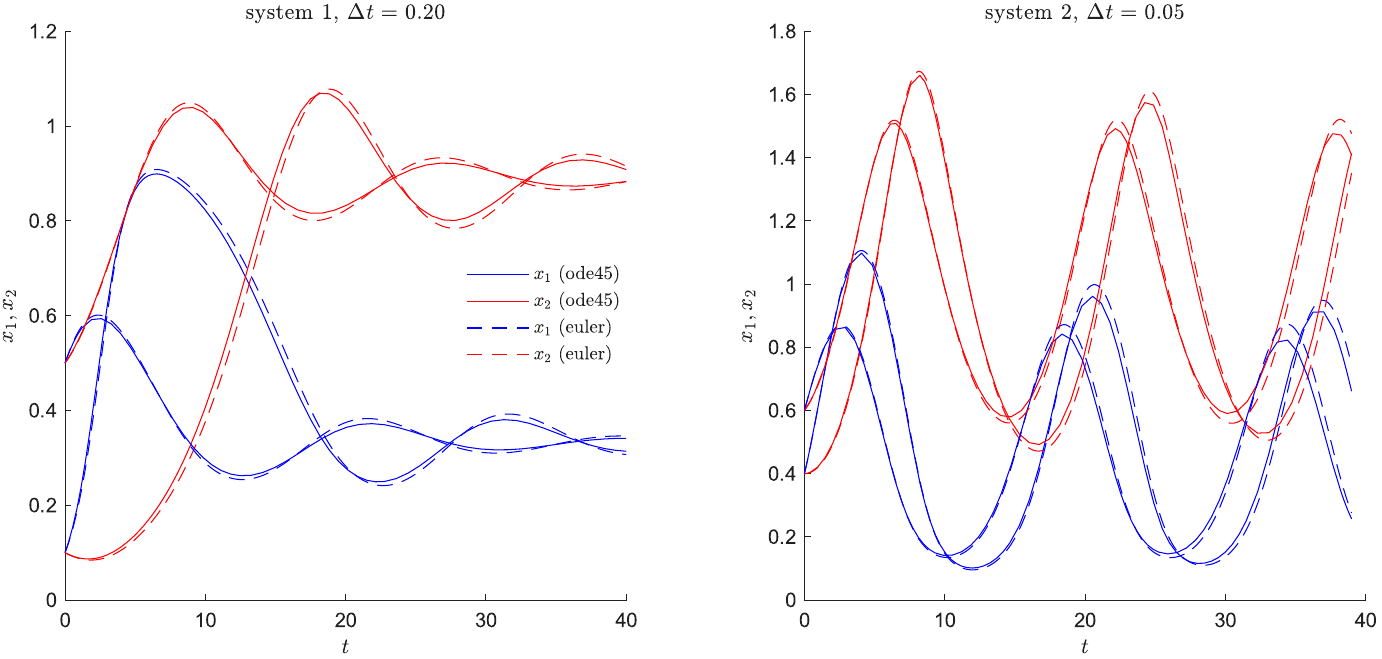}
\caption{Numerical impression of step sizes for the discretization of the dynamic systems by deterministic simulations in the relevant time and initial value range (two starting value examples, which are vertices of each initial density box) with \textit{ode45} (\textit{Runge-Kutta-(4,5)} algorithm) as ground truth and the \textit{Euler} algorithm which is basis of the simulation strategy. Left: for system 1 utilized in the simulation in Figure \ref{fig:Res:RosMcA_iteration_posterior_004_2} with starting values $(0.1,0.1)$ and $(0.5,0.5)$, right: for system 2 utilized in the simulation in Figure \ref{fig:Res:RosMcA_iteration_posteriorSDE_004} with starting values $(0.4,0.4)$ and $(0.6,0.6)$. Sufficient numerical approximation accuracy in the context of this study can be observed considering the simulated stochasticities. }  
\label{fig:Res:RosMcA_iteration_numericstep}
\end{figure}

In the second example, we present a SDE version of the normalized \textit{Rosenzweig McArthur} model \cite{Wang 2025}, i.e. 
\begin{align}
\text{d} x_1 & = \left(x_1\,\left(1-\frac{x_1}{k}\right)-\frac{m\,x_1\,x_2}{1+x_1}\right)\,\text{d} t + \sigma_1\,x_1\,\text{d}W_1\\
\text{d} x_2 & = \left(-c\,x_2 +\frac{m\,x_1\,x_2}{1+x_1}\right)\,\text{d} t + \sigma_2\,x_2\,\text{d}W_2\;,
\end{align}
with the real parameters $k,m,c$ and a Wiener process $\boldsymbol{W}_t$, leading to the transfer function 
\begin{align}
\boldsymbol{T}( \boldsymbol{x}^{(n)}, \boldsymbol{C}^{(n)} ) & := \left(\begin{array}{c} x_{1}^{(n)} +  \left( x_{1}^{(n)}\,\left(1-\frac{x_{1}^{(n)}}{k}\right)-\frac{m\,x_1^{(n)}\,x_2^{(n)}}{1+x_1^{(n)}} \right)\, \Delta t +  \sigma_1\,x_1^{(n)}\,C_1^{(n)}\\
x_{2}^{(n)} + \left( -c\,x_2^{(n)} +\frac{m\,x_1^{(n)}\,x_2^{(n)}}{1+x_1^{(n)}} \right) \, \Delta t +  \sigma_2\,x_2^{(n)}\,C_2^{(n)}\end{array}\right)\;.
\end{align}
In this simulation we utilize a parameter set which leads to periodical orbits in the $x_1$-$x_2$-space in the deterministic case with $k=1.9, m=1.1, c=0.31$, but add the stochastic noise densities $f_{C_1},f_{C_2} := \mathcal{N}(0,\Delta t)$ of the discretized SDE and use scaling coefficients $\sigma_1=\sigma_2=0.04$. {We utilize the same $f_{B_1},f_{B_2}$ as in the previous example.} In Figure \ref{fig:Res:RosMcA_iteration_posteriorSDE_004} we utilize $P=384000$ samples and $\Delta t=0.05$ (leading to the observation of a clear orbital structure with minor numerical artifacts, see Figure \ref{fig:Res:RosMcA_iteration_numericstep} (right) for an impression) and present in the plot representative density functions (iteration 476 to 780) roughly over one period of the orbit. The periodicity is not perfect since the parameter values are a compromise for the periodicity of the full starting value range (which lead also to slightly different periodicities, see Figure \ref{fig:Res:RosMcA_iteration_numericstep} (right)) and the SDE diffusion leads to broad distributions of the orbital state density. Still, a clear orbit is observable also in this SDE version and it is demonstrated that such challenging simulations can be performed based purely on the computational framework of iterative density transformations.

An additional verification of the full-density simulation is provided in Appendix Section \ref{sec:NumVec2DOrnUhl} where the simulation results are compared to the evolution of the analytically known mean and covariances of a 2D Ornstein-Uhlenbeck type SDE.

\subsection{Optimization with Full-Density Gradient Descent}
\label{res:FDGD}

\begin{figure}[htbp]
\centering
\includegraphics[width=14cm]{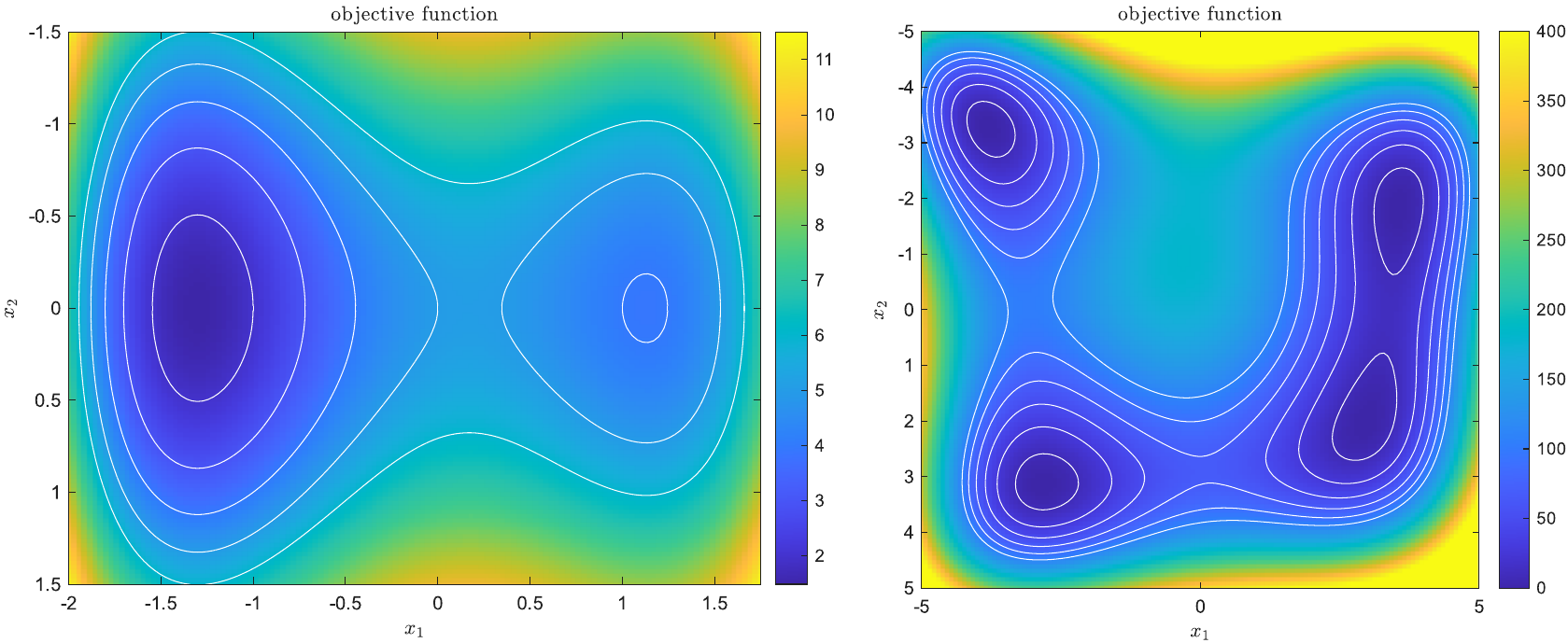}
\caption{Illustration of the objective function. Left: With two local minima on the $x_1$-axis. Right: \textit{Himmelblau's} function of optimization for four local minima.}  
\label{fig:Res:RGD_objective_001.pdf}
\end{figure}

\begin{figure}[htbp]
\centering
\includegraphics[width=14cm]{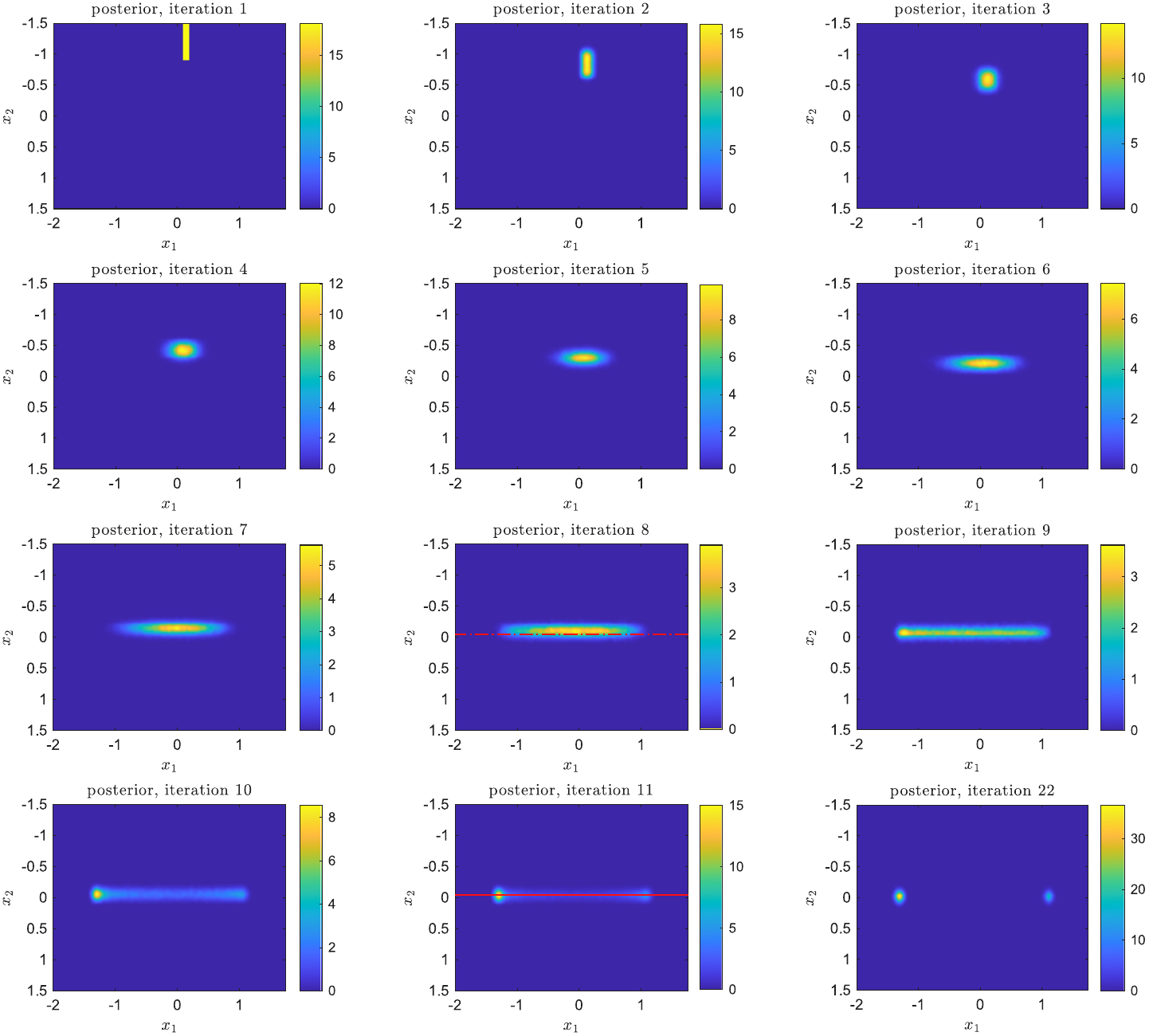}
\caption{Illustration of the posterior densities for different iteration numbers for FDGD-II for the objective function with two local minima. One can see the evolution of the square initial density being first attracted to the saddle point and second the convergence to the two local minima. {The computational density values were evaluated at $200\times 160$ grid points with equidistant spacing along each coordinate on $[-2,1.75]\times[-1.5,1.5]$. The initial density was uniform on the strip $[0.1,0.2]\times [-1.5,-0.9]$.}}  
\label{fig:Res:RGD_II_003}
\end{figure}

\begin{figure}[htbp]
\centering
\includegraphics[width=14cm]{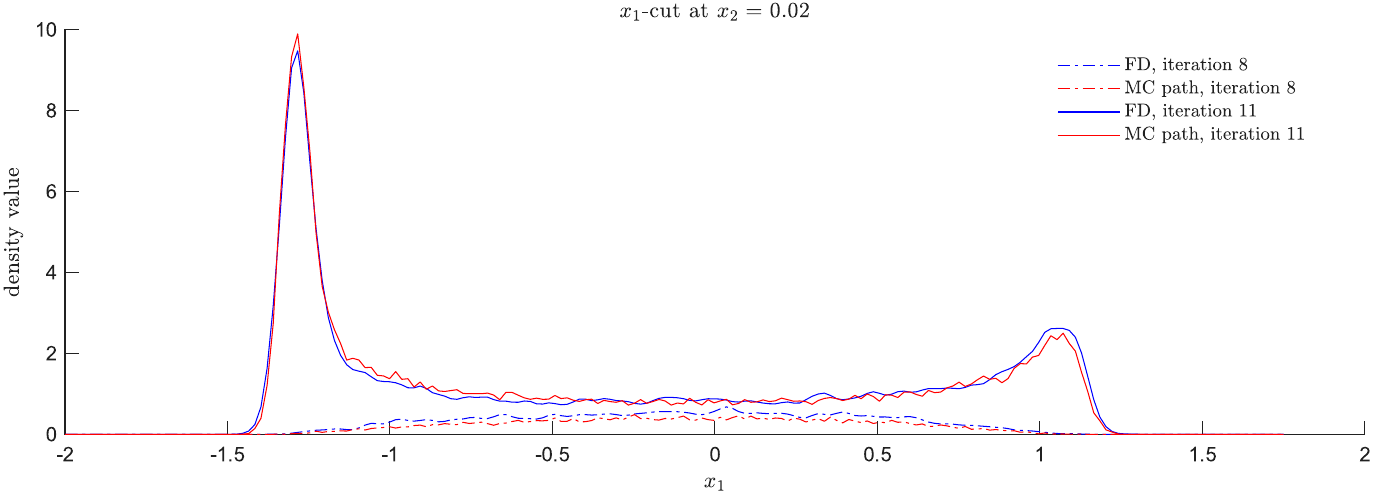}
\caption{{Verification of the computation of Figure \ref{fig:Res:RGD_II_003} by comparing cross sections of the numerical simulation of the full-density approach (FD) to the pathwise Monte Carlo simulation (MC path) utilizing $768000$ paths. The cross section are additionally presented in Figure \ref{fig:Res:RGD_II_003}. The two cross sections are $x_1$-cuts at $x_2=0.02$ for iterations $8$ and $11$.}  }
\label{fig:Res:RGD_II_003_cuts}
\end{figure}

\begin{figure}[htbp]
\centering
\includegraphics[width=14cm]{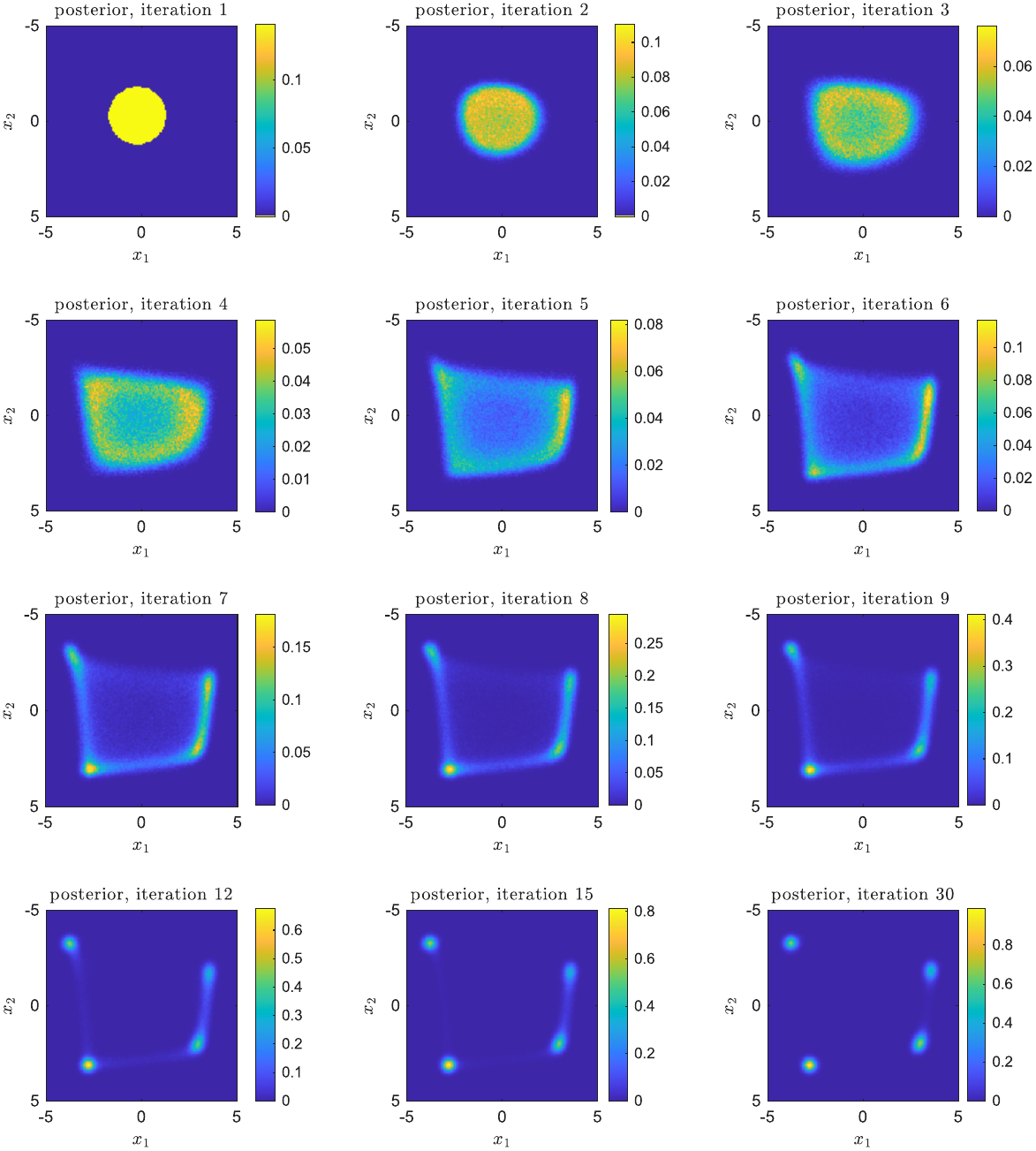}
\caption{Illustration of the posterior densities for different iteration numbers for FDGD-III using \textit{Himmelblau's} test objective function. One can see the evolution of the broad circular initial density with a variation of the update rate and its fast convergence to the four local minima. {The computational density values were evaluated at $200\times 200$ grid points with equidistant spacing on $[-5,5]\times[-5,5]$. The initial density was uniform on the circle with center $(-0.2,-0.3)$ and radius $1.5$.}}  
\label{fig:Res:RGD_III_012}
\end{figure}

{As the second application we are applying the presented FDGD optimization to two challenging types of objective functions in optimization: a saddle-point scenario and the \textit{Himmelblau's} function as a classical benchmark. This illustrates the use of the RIE modeling framework to gradient descent.} First, we present an objective function with two local minima and one saddle point in between
\begin{align}
F(\boldsymbol{x}) = x_1^4-3\,x_1^2+x_1+5+2\,x_2^2\;,
\end{align}
presented in Figure \ref{fig:Res:RGD_objective_001.pdf} on the left side with
\begin{align}
\nabla F(\boldsymbol{x}) = \left(\begin{array}{c} 4\,x_1^3 - 6\,x_1+1\\ 4\,x_2 \end{array}\right)\;,
\end{align}
and apply the FDGD-II algorithm, leading to the transfer function
\begin{align}
\boldsymbol{T}( \boldsymbol{x}^{(n)}, \boldsymbol{C}^{(n)} ) & := \left(\begin{array}{c} x_1^{(n)} - \eta\cdot\left(4\,\left(x_1^{(n)}\right)^3 - 6\,x_1^{(n)}+1  \right) + C_1^{(n)}\\ x_1^{(n)} - \eta\cdot\left(4\,x_2^{(n)}\right)  + C_2^{(n)} \end{array}\right)\;.
\end{align}
Second, we present the famously challenging \textit{Himmelblau's} objective function \cite{Himmelblau 1972, Przbylowicz 2026} with four minima which are difficult to localize
\begin{align}
F(\boldsymbol{x}) = (x_1^2 + x_2 - 11)^2 + (x_1+x_2^2-7)^2\;,
\end{align}
presented in Figure \ref{fig:Res:RGD_objective_001.pdf} on the right side with
\begin{align}
\nabla F(\boldsymbol{x}) = \left(\begin{array}{c} 4\,x_1\,(x_1^2 + x_2 - 11) + 2\,(x_1+x_2^2-7) \\ 2\,(x_1^2 + x_2 - 11) + 4\,x_2\,(x_1+x_2^2-7)  \end{array}\right)\;,
\end{align}
and apply the FDGD-III algorithm, leading to the transfer function
\begin{align}
\boldsymbol{T}( \boldsymbol{x}^{(n)}, \boldsymbol{C}^{(n)} ) & := \left(\begin{array}{c} x_1^{(n)} - C_3^{(n)}\cdot\left(4\,x_1^{(n)}\,\left(\left(x_1^{(n)}\right)^2 + x_2^{(n)} - 11\right) + 2\,\left(x_1^{(n)}+\left(x_2^{(n)}\right)^2-7\right)  \right) + C_1^{(n)}\\  x_1^{(n)} - C_3^{(n)}\cdot\left(2\,\left(\left(x_1^{(n)}\right)^2 + x_2^{(n)} - 11\right) + 4\,x_2^{(n)}\,\left(x_1^{(n)}+\left(x_2^{(n)}\right)^2-7\right)  \right)  + C_2^{(n)} \end{array}\right)\;.
\end{align}
In Figure \ref{fig:Res:RGD_II_003}, for the first objective function (Figure \ref{fig:Res:RGD_objective_001.pdf}, left) we utilize independent Gaussian densities $f_{C_1} := \mathcal{N}(0,0.04^2)$ and $f_{C_2} := \mathcal{N}(0,0.04^2)$  as well as $f_{B_1},f_{B_2} := \mathcal{N}(0,0.02^2)$ {(again truncated at three standard deviations)} and $P=48000$ samples.  This FDGD-II example with learning rate $\eta=0.075$ and a small starting density shows how the gradient descent is attracted by the saddle point first and afterwards converges simultaneously to the two minima. {In addition, in Figure \ref{fig:Res:RGD_II_003_cuts} cross sections of that simulation are presented and compared to pathwise Monte Carlo simulation. The shapes and positions of the intensity values are comparable and the smoothing effect of $\boldsymbol{B}$ is negligible.} In Figure \ref{fig:Res:RGD_III_012}, for a FDGD-III example with the second objective function (Figure \ref{fig:Res:RGD_objective_001.pdf}, right) we utilize $f_{C_1} := \mathcal{N}(0,0.2^2)$, $f_{C_2} := \mathcal{N}(0,0.2^2)$  as well as the same $f_{B_1},f_{B_2}$ as in the previous optimization example and the learning rate density $f_{C_3} := \mathcal{N}(0.01,0.003^2)$ (using $P=384000$ samples). Starting with a broad density, we see the fast attraction by a broad search for the posterior to the four minima. This shows, that in this stochastic regime a broad search strategy can be beneficial finding challenging extrema in a simultaneous fashion.

\subsection{Chaotic Random Iterations}
\label{sec:ChaotiCRIE}

\begin{figure}[htbp]
\centering
\includegraphics[width=14cm]{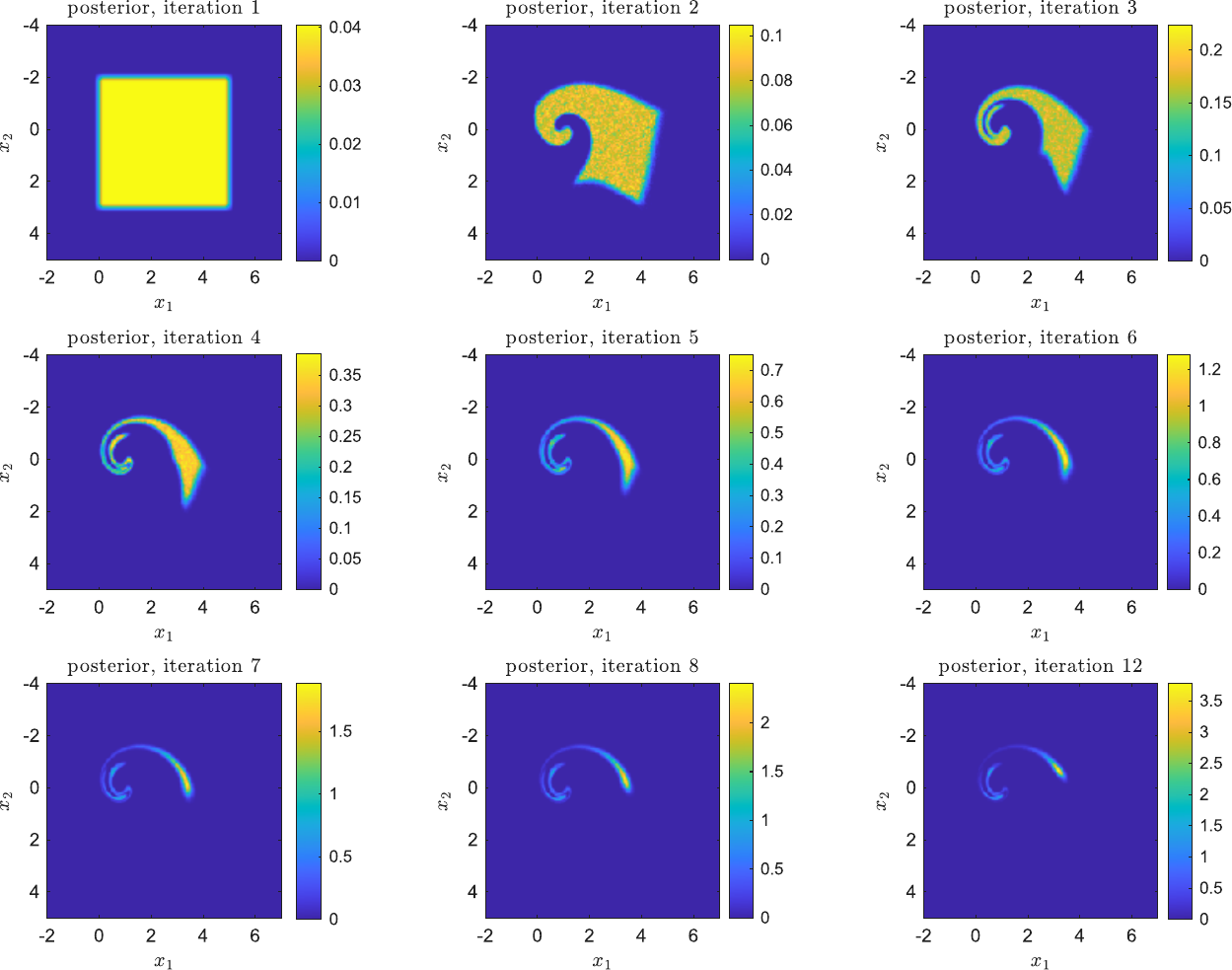}
\caption{Illustration of the posterior densities of the \textit{Ikeda mapping}. One can see the evolution of the broad square initial density into the characteristic strange attractor of the Ikeda mapping averaged over a density of the $u$-parameter. {The computational density values were evaluated at $200\times 200$ grid points with equidistant spacing on $[-2,7]\times[-4,5]$. The initial density was essentially uniform on the box $[0,5]\times[-2,3]$ with a slight Gaussian blur with standard deviation of two grid spacings.}}  
\label{fig:Res:StrangeAttract_001_2}
\end{figure}

{Finally, we also want to present how the RIE computational strategy can be applied directly to chaotic random iterations since this application class does not rely on a discretization procedure and the standard RIE formulation of Section \ref{sec:GenForm} can be directly utilized. We again present two numerical examples for this class of application.} First, the classical \textit{Ikeda map} \cite{Wang 2025 II} in real space is presented by
\begin{align}
x_1^{(n+1)} &= 1 + u\cdot ( x_1^{(n)}\cdot \cos(t^{(n)}) - x_2^{(n)}\cdot \sin(t^{(n)}))\\ 
x_2^{(n+1)} &= u\cdot ( x_1^{(n)}\cdot \sin(t^{(n)}) + x_2^{(n)}\cdot \cos(t^{(n)}))\;,
\end{align}
with
\begin{align}
t^{(n)} = 0.4 - \frac{6}{1+\left(x_1^{(n)}\right)^2+\left(x_2^{(n)}\right)^2}\;. 
\end{align}
Depending on the parameter $u$ different strange attractors appear. We are exploring the RIE for a small range of the variable $C_1^{(n)} = u$ represented by a random variable leading to the transfer function
\begin{align}
\boldsymbol{T}( \boldsymbol{x}^{(n)}, C_1^{(n)} ) & := \left(\begin{array}{c} 1 + C_1^{(n)}\cdot ( x_1^{(n)}\cdot \cos(t^{(n)}) - x_2^{(n)}\cdot \sin(t^{(n)}))  \\  C_1^{(n)}\cdot ( x_1^{(n)}\cdot \sin(t^{(n)}) + x_2^{(n)}\cdot \cos(t^{(n)})) \end{array}\right)\;.
\end{align}
In Figure \ref{fig:Res:StrangeAttract_001_2} we explore the progress of iterations for $f_{C_1} := \mathcal{N}(0.7,0.02^2)$ {as well as $f_{B_1},f_{B_2} := \mathcal{N}(0,0.02^2)$ (again truncated at three standard deviations)} with $P=384000$. It can be observed that the strange attractor arises during iterations, but also how the Ikeda mapping behaves in its first iterations acting on a full starting density distribution. This is an average of the strange attractor according to the ensemble distribution of the $u$-parameter (=$C_1^{(n)}$).

\begin{figure}[htbp]
\centering
\includegraphics[width=14cm]{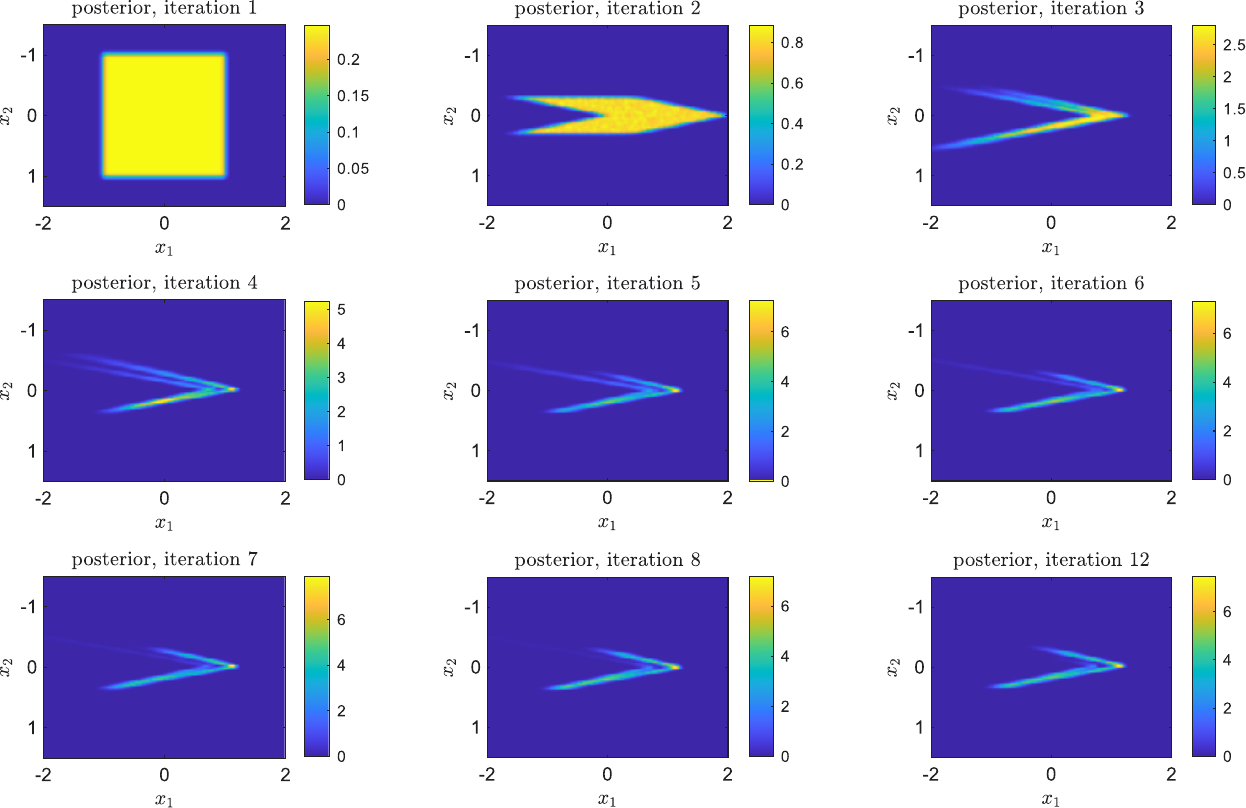}
\caption{Illustration of the posterior densities of the \textit{Lozi mapping}. One can see the evolution of the broad square initial density into the strange attractor of the Lozi mapping with sharps edges averaged over a density of the $a$-parameter. {The computational density values were evaluated at $200\times 150$ grid points with equidistant spacing on $[-2,2]\times[-1.5,1.5]$. The initial density was essentially uniform on the box $[-1,1]\times[-1,1]$ with a slight Gaussian blur with standard deviation of two grid spacings.}}  
\label{fig:Res:StrangeAttract_002}
\end{figure}

Second, the classical \textit{Lozi map} \cite{Lozi 2023} is presented, which is given by 
\begin{align}
x_1^{(n+1)} &= 1  - a\cdot |x_1^{(n)}| + x_2^{(n)}\\ 
x_2^{(n+1)} &= b\cdot x_1^{(n)}\;,
\end{align}
with typically a fixed standard parameter $b=0.3$ and a range for $a\in[1.4,1.7]$ for which strange attractors appear. This can be transferred to the transfer function 
\begin{align}
\boldsymbol{T}( \boldsymbol{x}^{(n)}, C_1^{(n)} ) & := \left(\begin{array}{c} 1  - C_1^{(n)}\cdot |x_1^{(n)}| + x_2^{(n)}  \\  b\cdot x_1^{(n)} \end{array}\right)\;,
\end{align}
being an example for a non-differentiable transfer function. In Figure \ref{fig:Res:StrangeAttract_002} we explore the progress of iterations for $f_{C_1} := \mathcal{N}(1.55,0.1^2)$ {with the same $f_{B_1},f_{B_2}$ as in the previous chaotic map} and $P=384000$. Also in this case the characteristic strange attractor with its sharp edges can be observed, even if we utilize a full ensemble of $a$ parameters.

\section{Discussion}
\label{sec:discussion}

{In this explorative study we investigated a modeling and computational framework for RIEs. On the one side, this work is based on the previous paper \cite{Hoegele 2026}, which demonstrated how to calculate the probability density of the solutions space of (static and nonlinear) random equation systems. On the other side, the main idea in this paper is to utilize this methodology in an iterative fashion, leading to dynamics of probability density propagations. This is a major extension and allows to model dynamical phenomena compared to static. Further, it is investigated how this new computational framework relates to standard representations such as the \textit{Perron-Frobenius} operator.} Essentially, the iteration transforms in each step the starting density function combined with the random variables contained in the iteration function itself in order to compute the next iteration state density. This can be regarded as computing \textit{ensemble solutions} for dynamical systems, e.g. \cite{Penenko 2025}. 

A major difference compared to pathwise Monte Carlo simulations is that the full density of the state space is calculated in each iteration step instead of several single paths which are averaged afterwards. This can be an advantage for the adaptive and more efficient sampling of the occurring densities during simulation, since in the pathwise approach only after the full simulations the densities appear statistically.

Applications to stochastic dynamical systems, such as RDEs and SDEs, demonstrate nontrivial evolutions of the state densities. We also show how this framework can be used for an explorative gradient descent, which we denote by FDGD (\textit{full-density gradient descent}) as a useful extension of classical gradient descent to \textit{global optimization}. Essentially, this leads to a meaningful search strategy for many local minima simultaneously with decisions taken not by single random paths but for a whole density of search directions. In this context, \textit{stochastic gradient descent} \cite{Goodfellow2016} could be associated with single random path of optimization and \textit{ensemble optimization} with the tracking of a number of sample paths in parallel during iterations, but not the full density. We want to emphasize that the presented simulation of dynamical systems and optimization are two main applications of applied mathematics, and utilizing RIEs allows a straightforward extension to their stochastic regimes. Finally, we provide an impression how strange attractors in chaotic mappings are computed with this framework for full densities of starting values and parameters. This is the most direct application of the framework.

Another important task in applied mathematics is solving nonlinear equation systems, e.g. by Newton's method or similar approaches, which typically lead in the deterministic case also to iteration equations. Since solving static equations (without a dynamic component) was already completely addressed in \cite{Hoegele 2026} by the general nonlinear form $\boldsymbol{M}(\boldsymbol{x};\boldsymbol{A}) = \boldsymbol{B}$ without the need of iteration equations, the extension to this computational framework is regarded as unnecessary and it is referred to the paper.

The most important point of the presentation of this framework is that it is simple (if there is a method to solve static random equations) and comparably general, and it could be applicable to many dynamics which are based on stochastic iterations (explicitly focusing on discretized continuous phenomena). {Possible larger future application areas could be: a) In dynamical systems the time step $\Delta t$ could also be taken as random variable (similar to the learning rate in FDGD-III), representing different discretization schemes of the dynamical system simultaneously. This leads to a new understanding of the simulation of dynamical systems, allowing to stochastically explore the discretization method itself. b) An important area concerns applications that combine stochastic dynamical systems and optimization, such as in \textit{stochastic control}. By computing full-densities in each update step of a control task it is possible to find the best actions that are considering the total uncertainty of the current state vector as well as the update equation. In general, combining the two areas, discretized dynamical systems and optimization within the RIE framework may lead to new insights in both areas.}

{
Of course, this study contains several crucial limitations which are briefly listed in this paragraph. a) The method is presented only for low dimensional examples. b) The presented accuracy of the examples depends on the discretization schemes (which we regard outside of the scope of the paper) as well as on the Monte Carlo integration (which is only partially verified, e.g. for the 2D Ornstein-Uhlenbeck SDE in Appendix \ref{sec:NumVec2DOrnUhl}, and by pathwise Monte Carlo simulations). c) This work does not provide formal convergence proofs for the general RIE computational strategy. d) The presented examples are illustrative and not exhaustive benchmarks. e) It is not tested how the framework behaves under stronger nonlinearities or higher-dimensional chaotic systems. f) The utilized acceptance rejection algorithm is only reasonably applicable in the low dimensional examples as presented. g) The choice of the regularization parameter, i.e. the standard deviation of the zero-centered random variable $\boldsymbol{B}$, is  manually chosen for each application. Smoothing effects are partially visible compared to pathwise Monte Carlo simulation and may play a role especially for large iteration numbers. In practice, a balance must be found between numerical diffusion due to the regularization and an adequate noise level as discussed in the numerical implementation Section \ref{sec:numimp}. h) Theoretical bounds or convergence rates of numerical schemes towards stationary densities are important on their own and considered complementary to this presentation. Instead this study explores the possibility to approximately compute and visualize the evolution of a finite number of iterations. We regard this focus as underrepresented in the study of stochastic dynamical systems which is a major part of the motivation of this work. Due to this specific focus, this approach should be regarded in the mindset of \textit{uncertainty quantification} in the broad field of \textit{scientific computing}, where uncertainties are propagated through dynamical systems by simulation and their impact is quantified in the results.
}

{Based on the limitations, there} are several dimensions to further develop the application of the framework: a) more efficient numerical schemes for sampling the density functions in the Monte Carlo integration than, e.g. \textit{acceptance rejection} sampling, b) utilizing the \textit{linear RIE operator} combined with a linear decomposition of the posterior density functions {such as presented in connection to Ulam's method in Appendix \ref{sec:AdvAnaRitz}}, c) more efficient use of the definition of the state space grid in the simulation {which were constant in this paper}, including possibly adaptive grids, d) better temporal discretization approaches in the derivation of the random iteration equations (e.g. for continuous dynamical systems from \textit{Euler} algorithm to \textit{Runge Kutta} algorithms, e.g. \cite{Amiri 2015}) and a detailed investigation of their convergence rates, e) extension of the optimization approaches, e.g. higher order iterative schemes such as full-density Newton optimization, etc., f) more challenging parameter density functions apart from normal distributions and a broader (also discontinuous) selection of the transfer functions, g) going from $R=2$ to higher dimensional state spaces and/or random variable vectors, e.g. investigating the \textit{curse of dimensionality} for this approach, and finally, h) extension to other interesting application areas where iteration equations can play a productive role in understanding dynamical phenomena and contain random variable densities as a new degree of freedom in the modeling process.

\section{Conclusion}

In this paper, a unified framework for the simulation of random iteration equations (RIE) via direct probability density propagation has been developed. This is based on recent work on static random equations, which is iteratively applied in this study. In consequence, the algorithmic approach is conceptually simple and straightforward. Expressive applications, such as the simulation of random and stochastic differential equations, the extension of gradient descent to the new full-density gradient descent approach as well as the simulation of chaotic maps with strange attractors showcase {possibly} wide applicability. {These examples illustrate complex stochastic evolutions and suggest that the framework may be useful for studying nonlinear dynamical systems under uncertainty, understanding the interplay between nonlinear dynamics and stochasticities. Several limitations of the presentation as well as future directions of this approach are discussed.}

\section{Appendix}

{
\subsection{Ritz-Galerkin Discretization and Ulam's Method for the RIE-Operator for Advanced Analysis}
\label{sec:AdvAnaRitz}

Starting from Equation \ref{equ:RIEIntOp}, we introduce the general Ritz-Ansatz for the density functions
\begin{align}
\pi_{\boldsymbol{x}^{(n)}}(\boldsymbol{x}) := \sum\limits_{i\in \Omega} \alpha_i^{(n)}\,\varphi_{i}(\boldsymbol{x})\;,
\end{align}
with a finite or infinite number of basis functions $\varphi_{i}$ of a function subspace in $L^2$, and the corresponding real coefficients $\alpha_i^{(n)}\in\mathbb{R}$. Please note, $\varphi_i$ are at this point arbitrary $L^2$ functions and depending on the numerical evaluation they need to be chosen in order to represent the density function meaningfully. Inserting this into Equation \ref{equ:RIEIntOp} leads to
\begin{align}
\pi_{\boldsymbol{x}^{(n+1)}}(\boldsymbol{x})\propto & \int\limits_{\mathbb{R}^R} \left(\; \int\limits_{\mathbb{R}^{K}}  f_{\boldsymbol{B}}(  \boldsymbol{x} - \boldsymbol{T}( \boldsymbol{s}_1, \boldsymbol{s}_2 ) ) \cdot f_{\boldsymbol{C}}(\boldsymbol{s}_2)\,\text{d}\boldsymbol{s}_2 \right)  \cdot\; \left( \sum\limits_{i\in \Omega} \alpha_i^{(n)}\,\varphi_{i}(\boldsymbol{s}_1) \right) \;\text{d}\boldsymbol{s}_1\\
= & \sum\limits_{i\in \Omega} \alpha_i^{(n)} \underbrace{ \left[\;\int\limits_{\mathbb{R}^R} \left(\; \int\limits_{\mathbb{R}^{K}}  f_{\boldsymbol{B}}(  \boldsymbol{x} - \boldsymbol{T}( \boldsymbol{s}_1, \boldsymbol{s}_2 ) ) \cdot f_{\boldsymbol{C}}(\boldsymbol{s}_2)\,\text{d}\boldsymbol{s}_2 \right)  \cdot\; \varphi_{i}(\boldsymbol{s}_1)  \;\text{d}\boldsymbol{s}_1\;\right]}_{=: p_i(\boldsymbol{x})} \;.
\end{align}
Utilizing the same Ritz-Ansatz for $\pi_{\boldsymbol{x}^{(n+1)}}(\boldsymbol{x})$, we get the brief expression
\begin{align}
 \sum\limits_{i\in \Omega} \alpha_i^{(n+1)}\,\varphi_{i}(\boldsymbol{x}) \propto \sum\limits_{i\in \Omega} \alpha_i^{(n)}\,p_{i}(\boldsymbol{x})\;.
\end{align}
Projecting these two linear decompositions onto the basis $\varphi_i$ itself utilizing the $L^2$ inner product (i.e. \textit{Galerkin}-type argument) leads for all $j\in\Omega$ to
\begin{align}
 \sum\limits_{i\in \Omega} \alpha_i^{(n+1)}\,\langle\, \varphi_{i}(\boldsymbol{x}),\varphi_{j}(\boldsymbol{x})\,\rangle \propto \sum\limits_{i\in \Omega} \alpha_i^{(n)}\, \langle\,p_{i}(\boldsymbol{x}),\varphi_{j}(\boldsymbol{x})\,\rangle \;,
\end{align}
or shorter if the number of basis functions is finite
\begin{align}
\left[\begin{array}{ccc} & \vdots & \\ \dots & \langle\, \varphi_{i}(\boldsymbol{x}),\varphi_{j}(\boldsymbol{x})\,\rangle & \dots \\ & \vdots & \end{array}\right]\,\boldsymbol{\alpha}^{(n+1)} \propto \left[\begin{array}{ccc} & \vdots & \\ \dots & \langle\, p_{i}(\boldsymbol{x}),\varphi_{j}(\boldsymbol{x})\,\rangle & \dots \\ & \vdots & \end{array}\right]\,\boldsymbol{\alpha}^{(n)}\;,
\end{align}
leading to finite dimensional quadratic matrices, which can be rearranged to 
\begin{align}
\boldsymbol{\alpha}^{(n+1)} \propto \underbrace{\left[\begin{array}{ccc} & \vdots & \\ \dots & \langle\, \varphi_{i}(\boldsymbol{x}),\varphi_{j}(\boldsymbol{x})\,\rangle & \dots \\ & \vdots & \end{array}\right]^{-1}\cdot \left[\begin{array}{ccc} & \vdots & \\ \dots & \langle\, p_{i}(\boldsymbol{x}),\varphi_{j}(\boldsymbol{x})\,\rangle & \dots \\ & \vdots & \end{array}\right]}_{=:P}\,\boldsymbol{\alpha}^{(n)}\;,
\end{align}
with the \textit{propagation matrix} $P$, transforming the coefficients of the density $\boldsymbol{\alpha}^{(n)}$ at iteration $n$ to $\boldsymbol{\alpha}^{(n+1)}$ at iteration $n+1$. This can be understood as the discretized stochastic \textit{Perron-Frobenius} operator utilizing the Galerkin approach. Practically, one would first calculate $P$ for a given RIE, then $P\cdot \boldsymbol{\alpha}^{(n)}$, and normalize the result, so that the density function based on $\boldsymbol{\alpha}^{(n+1)}$ is proper with probability mass $1$. The last step can be performed iteratively, which means, each RIE has a constant propagation matrix $P$, which contains all the information about the nonlinear dynamical system projected onto the basis $\varphi_i$. In consequence, spectral analysis, steady state investigation, efficient numerical computation and more can be performed by analyzing and/or approximating $P$ for meaningful basis functions $\varphi_i$.

For example, if one utilizes $\varphi_i(\boldsymbol{x}) := \delta(\boldsymbol{x}-\boldsymbol{x}_i)$, i.e. the basis function is a dirac distribution at position $\boldsymbol{x}_i$ as part of an approximation of a (fine, equidistant) spatial grid (where $i$ is the linear index running through that grid), then the propagation matrix itself is approximated by the entries
\begin{align}
P_{ji} =  \int\limits_{\mathbb{R}^{K}}  f_{\boldsymbol{B}}(  \boldsymbol{x}_j - \boldsymbol{T}( \boldsymbol{x}_i, \boldsymbol{s}_2 ) ) \cdot f_{\boldsymbol{C}}(\boldsymbol{s}_2)\,\text{d}\boldsymbol{s}_2 \;,
\end{align}
a straightforward integration, although typically for a large matrix $P$ depending on the number of grid points. Again, the role of $\boldsymbol{B}$ is a regularization that allows to calculate $P$ practically. With this derivation we essentially have arrived at a stochastic extension of \textit{Ulam's} method for investigating dynamic systems by the local transfer probabilities from one grid point to another, e.g. see \cite{Surasinghe 2024}.

}

\subsection{Numerical Verification with the 2D Ornstein-Uhlenbeck SDE}
\label{sec:NumVec2DOrnUhl}

\begin{figure}[htbp]
\centering
\includegraphics[width=14cm]{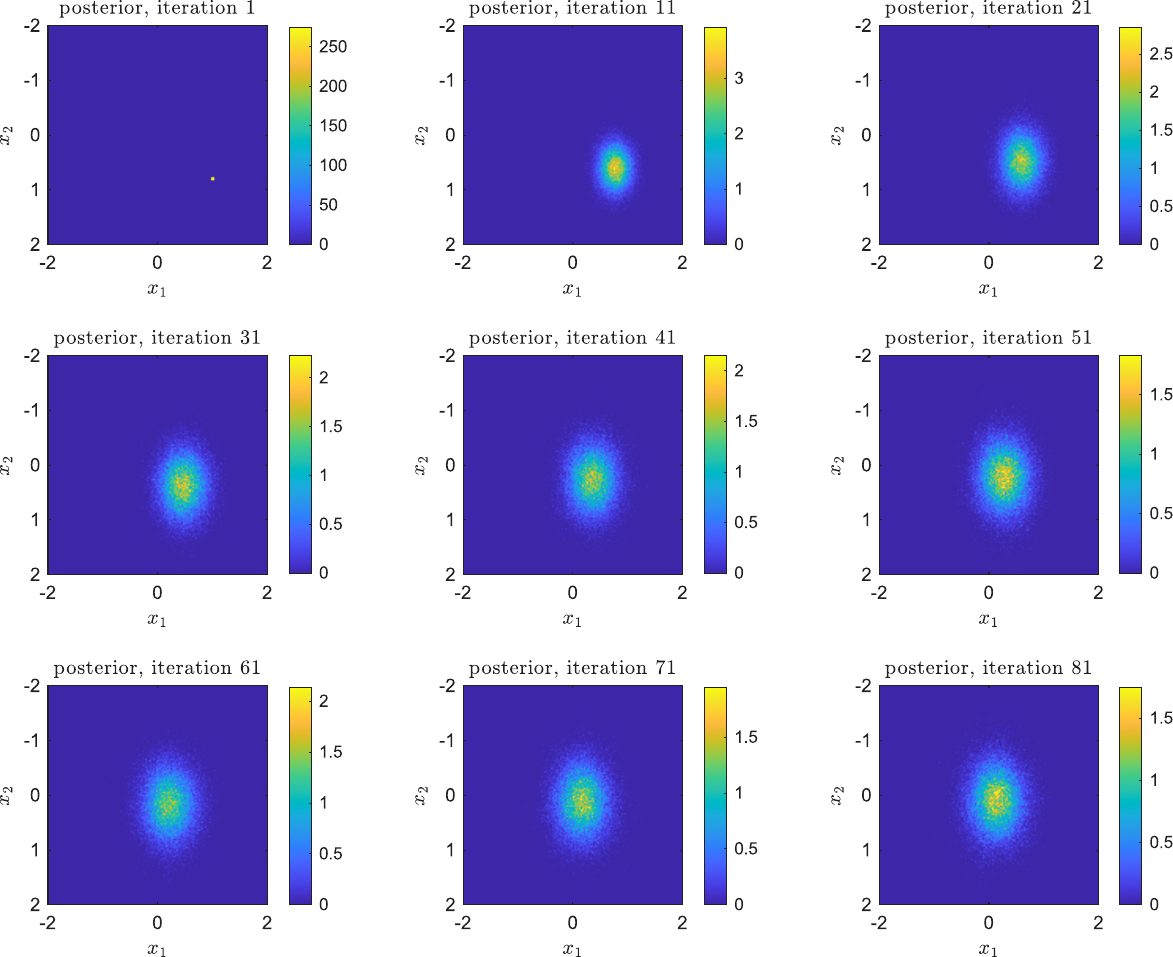}
\caption{Illustration of the posterior densities for different iteration numbers of the 2D Ornstein-Uhlenbeck SDE model utilizing the full-density approach. {The computational density values were evaluated at $200\times 200$ grid points with equidistant spacing on $[-2,2]\times[-2,2]$. The initial density was uniform on the box $[0.97,1.03]\times[0.77,0.83]$ (effectively a $2\times 2$ grid box).}}  
\label{fig:Res:2DOrnsteinUhlenbeck_verif_001}
\end{figure}

For an independent verification of the SDE calculation, we present a simple 2D Ornstein-Uhlenbeck SDE
\begin{align}
\text{d} x_1 & = -x_1\,\text{d} t + \sigma_1\,\text{d}W_1\\
\text{d} x_2 & = -x_2\,\text{d} t + \sigma_2\,\text{d}W_2\;,
\end{align}
with the real diffusion parameters $\sigma_1,\sigma_2$ and a Wiener process $\boldsymbol{W}_t$, leading to the transfer function 
\begin{align}
\boldsymbol{T}( \boldsymbol{x}^{(n)}, \boldsymbol{C}^{(n)} ) & := \left(\begin{array}{c} x_{1}^{(n)}    -x_{1}^{(n)} \, \Delta t +  \sigma_1\,C_1^{(n)}\\
x_{2}^{(n)}   -x_{2}^{(n)} \, \Delta t +  \sigma_2\,C_2^{(n)}\\ \end{array}\right)\;.
\end{align}
The stochastic noise densities are $f_{C_1},f_{C_2} := \mathcal{N}(0,\Delta t)$ of the discretized SDE and the diffusion coefficients are $\sigma_1=0.4,\sigma_2=0.6$. In Figure \ref{fig:Res:2DOrnsteinUhlenbeck_verif_001} we utilize $P=384000$ samples, $f_{B_1},f_{B_2} := \mathcal{N}(0,0.0025^2)$ {(truncated at three standard deviations)} and $\Delta t=0.025$ and show the evolution of the density function. As a starting distribution an (almost) point distribution at $\boldsymbol{x}_0 = (1,0.8)^T$ is chosen (please note, we utilized a slightly extended square with size $0.045$ in order to utilize exactly the same algorithm as in Section \ref{sec:res:dynamicsys}). During this transient phase, the expected drift to the zero mean position is observable as well as the anisotropic increase in the diffusion.\medskip

\begin{figure}[htbp]
\centering
\includegraphics[width=14cm]{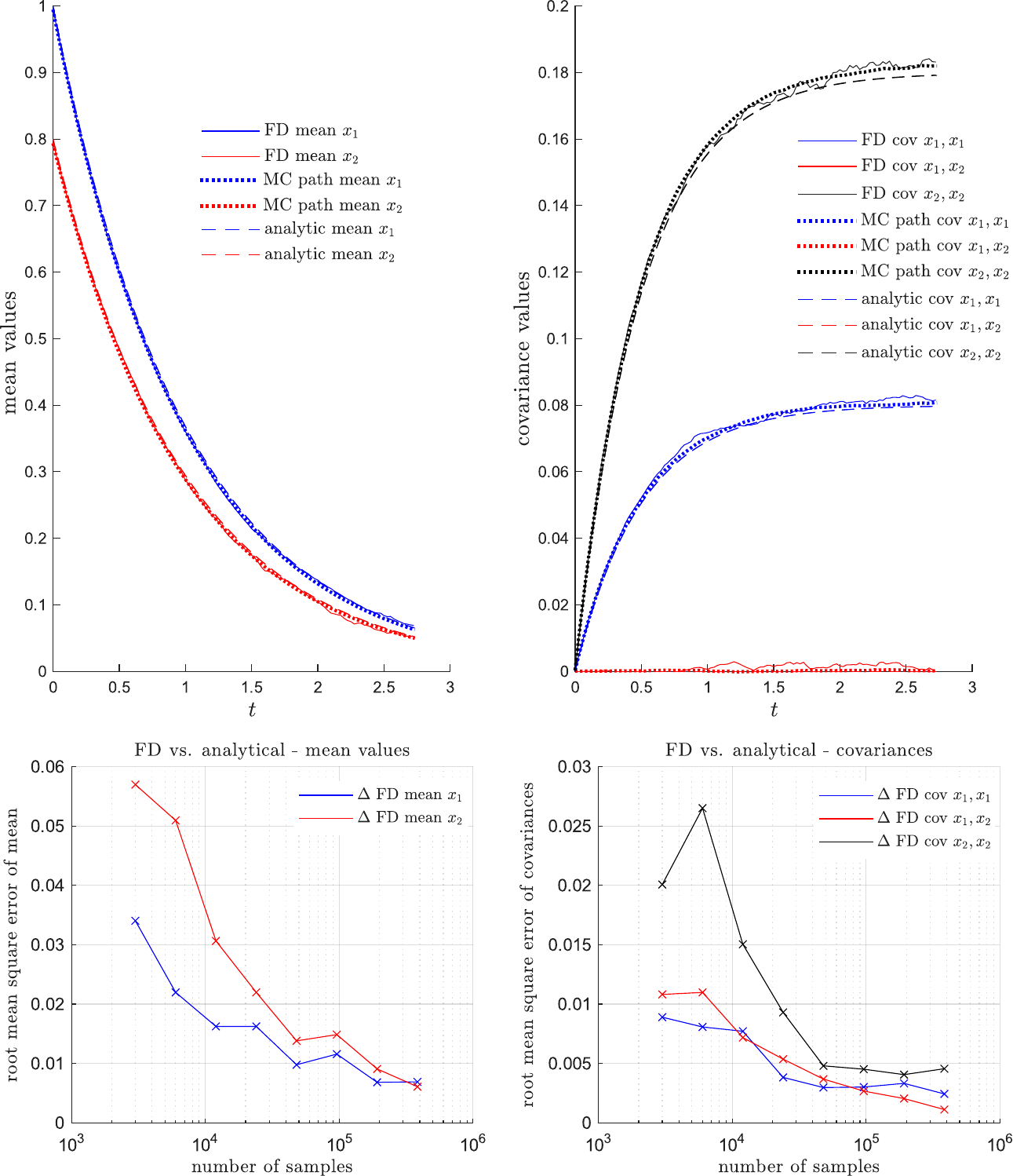}
\caption{Top row: Illustration of the mean and covariance values empirically evaluated for the full-density algorithm (FD, solid lines){, the pathwise Monte Carlo simulation (MC path, dotted lines)} and the analytic solution (dashed lines). The strong similarity between full-density (FD){, pathwise Monte Carlo (MC path)} and the analytic solution is demonstrated {verifying consistent calculations}. {Bottom row: Convergence analysis of the FD compared to the analytic solution is presented, showing the root mean squared error of the mean values (left) and covariances (right) in the same transient time frame $[0,2.725]$ as in the top row.}}  
\label{fig:Res:2DOrnsteinUhlenbeck_mean_cov_eval}
\end{figure}

The analytic solution for the mean $\boldsymbol{M}(t)$ and the covariance matrix $\boldsymbol{\Sigma}(t)$ for this process is known and given by
\begin{align}
\boldsymbol{M}(t) &= \boldsymbol{x}_0\,e^{-t}\\
\boldsymbol{\Sigma}(t) &= \frac{1}{2}\,(1-e^{-2\,t})\,\left(\begin{array}{cc} \sigma_1^2 & 0 \\ 0 & \sigma_2^2\end{array}\right)\;.
\end{align}
In Figure \ref{fig:Res:2DOrnsteinUhlenbeck_mean_cov_eval} \textit{top row} the numerical evaluation of the mean and the covariance values utilizing the RIE full-density approach over 109 iterations in the transient phase are compared to the time-dependent analytic solution {as well as to the pathwise Monte Carlo simulation with $384000$ single paths}. The comparison demonstrates an agreement with minor differences verifying the calculation. Possible sources for remaining errors are the discretization errors (spatial as well as temporal utilizing the Euler-Maruyama approach), the finite number of Monte Carlo samples for integration, an extended starting distribution in the simulation instead of a point distribution as for the analytic solution and the slightly diffusing character of the standard deviation of $\boldsymbol{B}$ when solving the random iteration equations. {It is observable that the influence of $\boldsymbol{B}$ is not dominant since pathwise Monte Carlo simulation (not having this regularization parameter) and the full-density approach agree with only minor fluctuations during this transient phase. In Figure \ref{fig:Res:2DOrnsteinUhlenbeck_mean_cov_eval} \textit{bottom row} a convergence analysis with respect to the number of Monte Carlo samples of $\boldsymbol{C}^{(n)}$ with $P\in\{3000,6000,12000,24000,48000,96000,192000,384000\}$ utilized in integration is presented. In these two figures the \textit{root mean squared errors} between the analytic solution and the full-density calculation of the mean values (left) and covariances (right) are presented in the time frame of $[0,2.725]$ (the same transient phase as in the top row) averaged over three runs. It is demonstrated how the simulation error reduces by increasing samples sizes suggesting consistency of the calculation. It must be noted, that minor differences due to the discretization scheme (see top row) will lead to a lower saturation bound.} 

\end{document}